\documentclass[11pt,fleqn]{article}
\usepackage{amssymb}
\usepackage{amsmath}
\usepackage{amstext}
\renewcommand{\thepage}{}
\renewcommand{\appendix}{\footnotesize\parindent .5cm\setcounter{equation}{0} 
\renewcommand{\theequation}{A.\arabic{equation}}
\setcounter{lemma}{0}\renewcommand{\thelemma}{A.\arabic{lemma}}}

\newtheorem{assumption}{Assumption}[section]
\newtheorem{hypothesis}{Hypothesis}

\newtheorem{proposition}{Proposition}[section]
\newtheorem{definition}{Definition}
\def\monthname{\ifcase\month\or
January\or February\or March\or April\or May\or June\or
July\or August\or September\or October\or November\or December\fi}

\newcommand{\pr}{{\rm pr}}
\newcommand{\elc}{{\rm elc}}
\newcommand{\corr}{{\rm score}}

\newcommand{\mme}{\mathbb{E}}
\newcommand{\mmx}{\mathbb{X}}
\newcommand{\mms}{\mathbb{S}}
\newcommand{\mmv}{\mathbb{V}}
\newcommand{\mmg}{\mathbb{G}}

\newcommand{\bondfowler}{{\rm B}}

\newcommand{\oy}{\overline{Y}}
\newcommand{\bh}{{\bf H}}
\newcommand{\bff}{{\bf F}}

\newcommand{\cale}{{\cal E}}

\newcommand{\obs}{{\rm obs}}
\newcommand{\direct}{{\rm direct}}
\newcommand{\spill}{{\rm spill}}

\newcommand{\bw}{{\bf W}}
\newcommand{\mmw}{\mathbb{W}}
\newcommand{\mmy}{\mathbb{Y}}
\newcommand{\mmr}{\mathbb{R}}
\newcommand{\mmc}{\mathbb{C}}
\newcommand{\mmp}{\mathbb{P}}
\newcommand{\bww}{{\bf w}}
\newcommand{\ba}{{\bf G}}
\newcommand{\bx}{{\bf X}}
\newcommand{\by}{{\bf Y}}

\begin{document}

\title{Exact P-values for Network Interference\thanks{\small We are grateful for comments by Peter Aronow, Peter Bickel, Bryan Graham, Brian Karrer, Johan Ugander, and seminar and conference participants at Cornell, the California Econometrics Conference, the UC Davis Institute for Social Science Inaugural Conference, and the Network Reading group at Berkeley.}
}
\author{Susan Athey\thanks{\small Graduate School of Business, Stanford University and NBER, athey@stanford.edu.
}
\and
Dean Eckles\thanks{\small Facebook, deaneckles@fb.com.
}
\and
 Guido W. Imbens\thanks{\small Graduate School of Business, Stanford University and NBER, imbens@stanford.edu.
}}
\date{
\monthname \ \number\year
}
\maketitle

\begin{abstract}
We study the calculation of exact p-values for a large class of non-sharp null hypotheses about treatment effects in a setting with  data from experiments involving members of a single connected network. The class includes null hypotheses that limit the effect of one unit's treatment status on another according to the distance between units; for example, the hypothesis might specify that the treatment status of immediate neighbors has no effect, or that units more than two edges away have no effect. We also consider hypotheses concerning the validity of sparsification of a network (for example based on the strength of ties) and hypotheses restricting heterogeneity in peer effects (so that, for example, only the number or fraction treated among neighboring units matters). Our general approach is to define an artificial experiment, such that the null hypothesis that was not sharp for the original experiment is sharp for the artificial experiment, and such that the randomization analysis for the artificial experiment is validated by the design of the original experiment.
\end{abstract}



\baselineskip=15pt\newpage \setcounter{page}{1}\renewcommand{\thepage}{[%
\arabic{page}]}\renewcommand{\theequation}{\arabic{section}.%
\arabic{equation}}

\section{Introduction}

This paper studies the calculation of exact p-values for a large class of non-sharp null hypotheses about treatment effects in a setting with  data from experiments involving members of a single connected network. For example, researchers might randomly assign some members, or clusters of members, of a social network to a treatment such as receiving information. 
We consider an environment where the following are observed: (i) the vector of treatments for all individuals in the network; (ii) the realized outcomes for all of the individuals (or possibly, a only a subset); (iii) all of the edges connecting individuals (where edges may potentially be categorized, for example into strong or weak edges); (iv) possibly, fixed characteristics for these individuals. Because the data come from a single network, with all units potentially connected and thus all units potentially affected by the full vector of treatments, establishing large sample approximations to distributions of statistics is challenging. This motivates our focus on the calculation of exact p-values based on the randomization distribution (Fisher, 1925). The validity of the p-value calculations does not depend on the network structure or the sample size.  Although we focus on the case of an explicit network where edges at most belong to a small
number of categories, the general methods we develop can be applied to more general settings with ``interference" and some measure of distance between units, where the researcher is interested in testing hypotheses about the nature of interference and how it relates to distance.

This paper considers a wide class of hypotheses about interference, sometimes caused by social interactions between units, where three categories of null hypotheses serve as leading examples.  In all three, the hypothesis
restricts the effects of the treatment of other units on a particular unit, while allowing for an individual's own treatment status to have a direct effect.  
The first category specifies that the treatment status of units with network distance weakly greater than \emph{k}
do not matter; when $k=1$, this requires that no other units' treatments have an impact, when $k=2$ only immediate neighbors' treatments matter, while when $k=3$ only neighbors as well neighbors of neighbors matter.   These types of hypotheses have been considered in empirical applications --- Bond et al. (2012) claim to find that ``messages not only influenced the users who received them, but also the users' friends, and friends of friends'' (p. 295) --- as well as in theoretical work, with many models a priori constraining spillovers in networks by ruling out effects of friends of friends (e.g., Toulis and Kao, 2013). The second category of hypotheses concerns the comparison between different categories of edges: e.g., under the null, only the treatments of neighbors with edges in one category matter. For example, Goldenberg, Zheng, Fienberg and Airoldi (2009) discuss a network defined through email interactions between Enron employees, with edges defined by the volume of email correspondence exceeding a threshold. 
Similarly, in analyses of large social networks researchers often sparsify the network by trimming edges between individuals with few interactions (see Thomas and Blitzstein, 2011, Bond et al., 2012, and Eckles, Karrer, and  Ugander, 2014). One can test whether such sparsification is appropriate by testing the hypothesis that  there are no spillovers between individuals not connected according to one definition of edges, but who would be connected under a looser definition of edges.  The third category of hypotheses concerns restrictions on heterogeneity
in the impact of neighbors.  For example, many models assume that only the number or fraction of treated neighbors matters for an individual's outcome, not which of their neighbors were treated.  An alternative of interest might be that neighbors with more connections matter more.   

There is a growing literature focusing on testing and inference in settings with general interference between units, both theoretical and empirical.\footnote{See Manski (1993, 2013), Christakis and Fowler (2007),  Rosenbaum (2007), Kolaczyk (2009), Aronow (2012),  
Bond, Fariss, Jones, Kramer, Marlow, and Fowler (2012), Bowers,  Fredrickson, and Panagopoulos (2012),  Hudgens and Halloran (2012),
Ugander, Karrer, Backstrom, and Kleinberg (2012),  Tchetgen and
 VanderWeele (2012), Goldsmith-Pinkham and Imbens (2013),  Liu and Hudgens (2013), Aronow and Samii (2014), Choi (2014), Eckles, Karrer, and Ugander (2014), Ogburn and VanderWeele (2014) and van der Laan (2014).} However, there is no available general asymptotic theory that handles hypothesis tests about these categories of null hypotheses, and
the nascent literature on estimation in network settings
requires strong restrictions on the network size and structure.\footnote{A small literature has emerged
that posits a specific functional form model of network formation (and thus the process for how the network changes as the size of the network grows), and then proposes an approach
for estimating the parameters of the network formation process (as opposed to parameters describing treatment effects).  In a leading example, Chadraskhar and Jackson (2014) establish consistency and asymptotic normality 
of parameter estimates for network formation under certain conditions (e.g. network is sufficiently sparse for a class of models they call subgraph generation models).  See also 
Holland and Leinhard (1981), Kolaczyk (2009), Manski (1993, 2013), Goldsmith-Pinkham and Imbens (2013), and Aronow and Samii (2014).}  In contrast,
our primary goal is to test hypotheses about the impact of treatments in a network setting, without restricting the network. 
  
The main contribution of this paper is to expand the applicability of the ``randomization inference'' approach to calculating exact p-values, originally developed by Fisher (1925) and Rosenbaum (1984), to our hypotheses of interest.\footnote{For applications of randomization inference outside the network setting, see Basu (1980), Rubin (1980), Rosenbaum (1995, 2002, 2007,  2009), Lehmann and Romano (2005), Imbens and Rubin (2015), and Canay, Romano, and Shaikh (2015).} In the randomization inference approach, the distribution of a test statistic is generated by the assignment mechanism, keeping fixed the potential outcomes and characteristics of the units.  This approach only applies directly to ``sharp'' null hypotheses, whereby the null hypothesis allows the analyst to infer the outcomes of individuals under 
alternative (counterfactual) treatment vectors.  For example, the null hypothesis that the treatment has no effect whatsoever is sharp, because an individual's outcome is known (and equal to his
realized outcome) under all possible treatment vectors.  Given this, it is possible to simulate draws from the random assignment of treatment vectors, and calculate the test
statistic of interest under each simulated draw (in this example, a natural test statistic is the average difference in outcomes between treated and control individuals).  The distribution of these simulated test statistics converges to the true distribution of the test statistic as the number of draws grows, and this true distribution is exact for the given network size and
structure rather than a large sample approximation.  Thus, exact p-values for the null hypothesis of no treatment effects
can be derived in a network setting using a conventional application of randomization inference.

In contrast, the three leading categories of null hypotheses outlined above are not sharp because under the null hypotheses we cannot infer the exact values for all outcomes for all possible values of the treatment; since all of the categories allow the treatment to have a direct effect on individuals, their outcomes cannot be inferred under alternative treatment
assignments.  In this paper we present a novel approach to dealing with such non-sharp null hypotheses.  Closest in spirit to this paper, Aronow (2012) adapts the 
randomization inference approach to consider the specific
non-sharp null hypothesis that only an individual's own treatment and that of his immediate neighbors matter, corresponding to the first category described above with $k=1$.  Here, we
provide a general framework that applies to a much larger class of non-sharp null hypotheses.

At an abstract level we address the problem that the null hypothesis of interest is not sharp by introducing the notion of an artificial experiment that  differs from the experiment that was actually carried out. This artificial experiment will be chosen so that the randomization analysis we propose for it is validated by the design of the experiment that was actually carried out, and at the same time the null hypothesis of interest that was not sharp for the actual experiment, is sharp for the artificial experiment.
In simple settings this idea of analyzing an experiment that differs from the experiment that was actually carried out is often used implicitly. Suppose we have an experiment where for each unit in a finite  population a coin is flipped to determine the treatment assignment for that unit. Given the data, we may analyze the data as if the total number of treated units is fixed, whereas in the actual experiment the number of treated  units is random. Analyzing the experiment as if the number of treated units is fixed is valid because we can think of the original experiment being a sequential one where in the first stage the number of treated units is determined by a sequence of coin tosses, and in the second stage the the fixed number of treated units is selected from the population at random. The artificial experiment is now simply the second stage of the original experiment, conditional on the first stage. In this case there is no loss of information because the number of treated units is ancillary.\footnote{This is similar in spirit to Rosenbaum (1984), who carries out randomization tests conditional on covariates or functions thereof such as the propensity score.}

In the settings we analyze in the current paper we also decompose the original experiment into two stages, and we analyze the experiment performed in the second stage conditional on the first stage randomization.  In an additional modification to the original experiment, we focus on a limited set of test statistics, namely those that depend on outcomes only for a subset of the original population, which we call the ``focal units.'' These changes to the original experiment lead to an artificial experiment where the null hypothesis that is not sharp in the original experiment is sharp for the artificial experiment, and where randomization inference is validated by the original experiment. The choice of the focal units on whose outcomes the statistic may depend and the decomposition of the original experiment into two stages are intricately linked to  achieve the goal of defining an experiment with a sharp null hypothesis amenable to randomization inference. 

The choice of focal units will matter for the power of the tests, but for any choice of focal units our approach will lead to exact p-values. Given the choice of focal units, we derive the unique partition of the space of assignments into subsets such that the null hypothesis implies that the outcomes for all focal units must be constant within these subsets. Then the original experiment is re-interpreted as a sequential experiment where in the first stage the subset into which the assignment falls is determined, and in the second stage the assignment is drawn randomly from within the subset (with the likelihood of each assignment implied by the original experiment). The analysis of our artificial experiment then focuses on a test statistic constructed from outcomes for the subpopulation of focal units and relies on the second stage randomization, conditional on the randomization in the first stage, to construct the p-value for the test statistic.

With our framework for testing in place,
it is then possible to compare the statistical power of alternative test statistics.  We do this for our three categories of hypotheses,
and we propose statistics that will be optimal for particular models of network interactions.  This in turn lays the groundwork for future research about optimal experimental design when the goal
is to test a given hypothesis or set of hypotheses.

The remainder of the paper is organized as follows. In the next section we introduce the general set up and notation. In Section \ref{hypotheses} we discuss a number of the hypotheses that we consider. This is not an exhaustive list, but it contains what we view as leading examples of the hypotheses researchers may wish to consider in network settings. Section \ref{artificial} contains a general discussion of the notion of artificial experiments that lies at the heart of our approach.
In the next four sections, Sections \ref{nospillover}, \ref{higher}, \ref{sparse}, and \ref{heterogeneity} we discuss in detail how the approach would be implemented for the main categories of null hypotheses we consider. These details include discussions of the decisions researchers need to make regarding the choice of focal units and test statistics. In Section \ref{simulations} we present the results from some simulations to evaluate the statistical power of the tests for alternative statistics. Section \ref{conclusion} concludes.

\section{Set Up}
\label{section:setup}

We have information on a population $\mmp$ of $N$ individuals, with $i$ indexing the individuals. 
We also have a set of treatments $\mmw$. In most of our examples each individual is either exposed to an intervention of not, although that is not necessary for some of the results. In that case for individual $i$ the exposure is denoted by $W_i\in\{0,1\}$, with $\bw$ the $N$-component vector of exposures with $i$th component equal to $W_i$, and $\mmw=\{0,1\}^N$.
There is mapping $\by:\mmw\mapsto \mmy^N$ of potential outcomes, with the $i$th element of this mapping written as
$Y_i:\mmw\mapsto \mmy$, where  $\mmy\subset \mmr$ is the set of values for the potential outcomes. We refer to $Y_i(\bww)$ as a potential outcome, with the corresponding vector of potential outcomes denoted by $\by(\bww)$. 
For the realized value of the assignment  $\bw$ we observe the corresponding vector of potential outcomes,
\[ \by^\obs=\by(\bw).\]
The treatment exposure $\bw$ is assigned through  an assignment mechanism $p:\mmw\mapsto [0,1]$,
where $p(\bww)$ is the probability of the assignment $\bw$ taking  on the value $\bww$, $p(\bww)={\rm pr}(\bw=\bww)$, satisfying $p(\bww)\geq0$ and $\sum_{\bww\in\mmw} p(\bww)=1$.

The units are connected through a undirected network that is observed by the researcher. The symmetric $N\times N$ adjacency
matrix $\ba$ measures the network, with the $(i,j)$th element of the adjacency matrix, denoted by $G_{ij}$, equal to one if there is an edge between units $i$ and $j$, and zero otherwise. By convention all diagonal elements  $G_{ii}$ are equal to zero.
We will call individuals $i$ and $j$ neighbors or peers if $G_{ij}=1$. 
The network is taken here to be a fixed characteristic of the population.
Let the distance $d(i,j)$ between units $i$ and $j$ be length of the shortest path between $i$ and $j$, and equal to $\infty$ if there is no path between $i$ and $j$. Thus $d(i,i)=0$, and $d(i,j)=1$
if $i\neq j$ and $G_{ij}=1$, $d(i,j)=2$ if $G_{ij}=0$ if $i\neq j$ but there is at least one unit $k$ such that $G_{ik}=1$ and $G_{kj}=1$, { et cetera}.
A special case is that with non-overlapping peer groups, considered, for example, in Manski (1993, 2013), Hudgens and Halloran (2008), and Carrell, Sacerdote and West (2013), where for all triples $(i,j,k)$,
$G_{ij}=1$ and $G_{jk}=1$  implies $G_{ik}=1$. We allow for such settings, but do not require them. Let $\mmg$ be the space of possible adjacency matrices.

For each individual there is also a vector of attributes $X_i$, with the  matrix of attributes denoted by $\bx$. Both the network and the attributes are viewed as pretreatment variables, not affected by the treatment.
We focus on the case where we observe the quadruple $(\by^\obs,\bw,\ba,\bx)$. More generally we may observe outcomes for a subset of the population.
The first two components of this quadruple, $\by^\obs$ and $\bw$ are random because of the randomization, the last two, $\ba$ and $\bx$, as well as the potential outcome function $\by(\cdot)$ are taken as fixed.

Let us think of an experiment for causal effects, denoted by $\cale$, being defined by a combination of the set $\mmw$  of possible values for the treatment $\bw$; the  population $\mmp$ of units characterized by their potential outcomes, their network and their fixed attributes; and a distribution for the treatment assignment, $p:\mmw\mapsto [0,1]$,
so that $\cale=(\mmw,\mmp,p(\cdot))$.

\section{Hypotheses}
\label{hypotheses}

In this section we discuss the three general classes of null hypotheses we consider, as well as some specific examples, and briefly discuss how p-values are calculated given a sharp null hypotheses. The classes  of hypotheses are not exhaustive, but they include many of the hypotheses that we view as interesting in settings with networks and are suggestive of the generality of the approach.

\subsection{Some General Concepts}
\label{general}

Let us start by formally defining several concepts: (i) a null hypothesis on treatment effects; (ii) whether a null hypothesis on treatment effects is \emph{sharp}; (iii) \emph{level sets}, that is, sets of assignments that result in invariant outcomes for a given individual.

\begin{definition}{\sc (A Null Hypothesis on Treatment Effects)}
A null hypothesis on treatment effects $H_0$  is a set of restrictions on the potential outcome function $\by:\mmw\mapsto \mmy^N$.\end{definition}
These restrictions can include the absence of any treatment effects, e.g., $Y_i(\bww)=Y_i(\bww')$ for all $\bww$, $\bww'$ and all $i$. They can also include more limited restrictions on the potential outcome functions.

\begin{definition}{\sc (A Sharp Null Hypothesis on Treatment Effects)}
A null hypothesis on treatment effects $H_0$ is sharp for $(\mmw,\mmp)$  if, given the value of  $(\bww,\by(\bww))$ for a single  assignment $\bww\in\mmw$, under $H_0$ we can infer 
 the value of $\by(\bww')$ for any other $\bww'\in\mmw$.\end{definition}

Now consider a test statistic, $T: \mmy^N\times \mmw\mapsto \mmr$. For a given experiment $\cale=(\mmw,\mmp,p(\cdot))$ the test statistic $T(\by(\bw),\bw)$ is random only through its dependence on the treatment (directly, and indirectly through the dependence of the realized outcome on the treatment). We can infer the distribution of the test statistic for a sharp null hypothesis. The p-value for the statistic under the null hypothesis is then the probability that the realization of the test statistic is at least as extreme as the observed value:
\[ \text{p-value}={\rm pr}\Bigl(|T(\by(\bw),\bw)|\geq |T(\by(\bw^\obs),\bw^\obs)|\Bigr).\]
In most cases we do not have available a closed form expression to calculate this p-value exactly. However, we can
 approximate it arbitrarily accurately by taking $B$ independent draws $\bw_b$ from the distribution of the assignment, $p(\cdot)$, and calculating the proportion of these $B$ draws that would have led to value for the statistic larger than or equal to the observed value of the statistic:
\[ \widehat{\text{p-value}}=\frac{1}{B}\sum_{b=1}^B\Bigl(|T(\bw_b,\by(\bw_b))|\geq |T(\bw^\obs,\by(\bw^\obs))|\Bigr),\]
for some large value of $B$. This estimate is unbiased for the true p-value, and its variance is bounded by $1/(4B)$, which can be made arbitrarily small by choosing $B$ large enough.

In some cases the statistic does not have a symmetric distribution under the null, and we may look at twice the minimum of the tail probabilities,
\[ \widehat{\text{p-value}}=
2\times \min\Big\{{\rm pr}\Bigl(T(\by(\bw),\bw)\geq T(\by(\bw^\obs),\bw^\obs)\Bigr)
,\]
\[\hskip3cm {\rm pr}\Bigl(T(\by(\bw),\bw)\leq T(\by(\bw^\obs),\bw^\obs)\Bigr)\Bigr\}
.\]

Most of the null hypotheses we consider in this paper are not sharp. However, they  imply that only a limited set of changes in the treatment actually change outcomes. To capture this, it is useful to introduce the notion of level sets, that is, sets of assignments with zero treatment effects.
\begin{definition}{\sc (Level Sets)}
Given a null hypothesis $H_0$, for each individual $i$  and for each treatment level $\bww$, define the level set $\mmv(i,\bww,H_0)$ as follows:
\[ \mmv(i,\bww,H_0)=\{\bww'\in\mmw| Y_i(\bww')=Y_i(\bww)\ {\rm given\ } H_0\}.\]

\end{definition}
Thus, the level set for unit $i$ given treatment vector $\bww$ is the set of treatments $\bww'$ such that under the null hypothesis, the potential outcome for unit $i$ is identical to the potential outcome given treatment $\bww$.\footnote{Manski and Tamer (2002) make use of level sets for non-network data. Related work on networks makes use of some concepts directly related to level sets. Manski (2013) and Eckles, Karrer, and Ugander (2014) work with \emph{effective treatments}, where each effective treatment corresponds to a level set, one of which is the observed level set. Aronow and Samii (2013) and Ugander et al. (2013) work with \emph{exposure models}, which uniquely specify effective treatments.}
(More generally we could define this set as the set of treatments where we can infer the potential outcomes, but outside the case where these potential outcomes are equal there are few cases of interest so we do not include that level of generality.)

These level sets play an important role in our approach, and it is useful to see what form they can take.
For the sharp null hypothesis that there is no treatment effect whatsoever, $\mmv(i,\bww,H_0)$ is equal to $\mmw$ for all $i$ and all $\bww$. With non-sharp null hypotheses, however, the set $\mmv(i,\bww,H_0)$ may vary, both  by individual and by treatment.  For example, in the setting where $\mmw=\{0,1\}^N$, if the null hypothesis allows for a direct effect of an individual's own treatment, but not for any effects of other individuals' treatment status, the set $\mmv(i,\bww,H_0)$ equals $\{\bww'\in\mmw|\bww'_i=\bww_i\}$, so that for each individual there are two possible values for the set, depending on the individual's own treatment status. At the other extreme, if the null hypothesis does not impose any restrictions, then level sets consist of singletons: $\mmv(i,\bww,H_0)=\{\bww\}$.
Because within a level set the treatment effect is zero, we can in principle do randomization inference on treatment effects for that individual.

It will play an important role later that in general for each unit $i$ these level sets define a partition of the assignment space $\mmw$ into $J$ level sets $\mmw_1,\ldots,\mmw_J$ such that for all $\bww\in\mmw$, $\mmv(i,\bww,H_0)\in\{\mmw_1,\ldots,\mmw_J\}$. If there are no restrictions at all, the elements of this partition consist of singletons, but in many interesting cases the number of elements of this partition is small. For example, for the null hypothesis that there are no spillovers, the partition contains two sets.

\subsection{Null Hypotheses on Spillovers}

We are interested in testing for the effect of exposure to the treatment for some individuals on the outcomes for others. We refer to such effects as ``spillovers,'' ``interactions'' or ``peer effects."  In the case where they are limited to the effects of direct neighbors, the peer effects we study are what Manski (1993) calls ``exogenous peer effects.''

First we consider the following three specific hypotheses that allow for a range of spillovers. Recall that in general the hypotheses we consider are restrictions on the mapping $\by:\mmw\mapsto \mmy^N$. 
\begin{hypothesis}\label{noeffects}{\sc (No Treatment Effects)}
 $Y_i(\bww)=Y_i(\bww')$ for all  $i$, and all pairs of assignments $\bww,\bww'\in\mmw$.
\end{hypothesis}
This is a sharp null hypothesis in the original experiment, because for all $\bww'\in\mmw$ the potential outcomes $Y_i(\bww')$ can be inferred from the observed treatment and observed outcomes $(\bww,\by(\bww))$ under the null hypothesis. Thus, the calculation of Fisher exact p-values is conceptually straightforward. 

Next, we consider a weaker null hypothesis that allows for effects of the own treatment on the own outcome, but not of the own treatment on a neighbor's outcome:
\begin{hypothesis}\label{nopeer}{\sc (No Spillovers)}
$Y_i(\bww)=Y_i(\bww')$ for all $i$, and all pairs of assignment vectors $\bww,\bww'\in\mmw$ such that
$\bww_i=\bww_i'$. 
\end{hypothesis}
This  null hypothesis is the one considered by Aronow (2012). It is not sharp, because it does not rule out that exposure to the treatment affects the outcome for the unit exposed. Manski (2013) refers to settings where this hypothesis holds as settings with ``individualistic treatment response.''
This null hypothesis is implied by the stable-unit-treatment-value-assumption (SUTVA, Rubin, 1980). Under this assumption we can simplify the notation to the conventional one in the causal effect literature where the potential functions are a function of the own treatment only, $Y_i(\bww)=Y_i(\bww_i)$. Because we consider more general cases, we continue to write the potential outcomes as a function of the full $N$-component vector $\bww$.


We can go beyond hypotheses ruling out all spillover effects, and allow for first order, but not higher order, spillover effects. That is, changing the treatment for neighbors may affect one's outcome, but changing the treatment for neighbors-of-neighbors does not change one's outcome.
\begin{hypothesis}\label{nohigher}{\sc (No Second and  Higher-Order Spillovers)}
 $Y_i(\bww)=Y_i(\bww')$ for 
all  $i$, and for all pairs of assignment vectors $\bww,\bww'\in\mmw$ such that
$\bww_j=\bww'_j$ for all units $j$ such that $d(i,j)<2$. 
\end{hypothesis}

Consider the following example where testing for higher order spillovers may be interesting.  
Suppose one can observe one's own treatment as well as the treatment of one's network neighbors, for example because of face-to-face interactions. One can also observe one's own outcome, but not the outcome for neighbors.
It may well be that in such cases there are spillover effects from neighbors, but no spillover effects from neighbors-of-neighbors or individuals even more distant in the network.
Testing for higher order spillover effects could then be interpreted as testing whether the network captures all the connections.

Some theoretical models (e.g. Toulis and Kao, 2013) model spillover effects in way that rules out higher order spillover effects. At the same time some researchers claim to find higher order spillovers effects in empirical work (e.g., Bond et al., 2012). Our tests are the first exact tests available for such hypotheses.


We can embed  these three hypotheses in a more general one that restricts $k$-th order spillover effects for arbitrary $k$. 
\begin{hypothesis}\label{kthorder}{\sc (No $k$-th and Higher Order Spillovers)}
For  unit $i$, for $i=1,\ldots,N$, $Y_i(\bww)=Y_i(\bww')$ for 
 all pairs of assignment vectors $\bww,\bww'\in\mmw$ such that
$w_j=w'_j$ for all units $j$ such that $d(i,j)<k$. 
\end{hypothesis}
(Here we interpret the set of pairs $\bww$ and $\bww'$ such that $\bww_i=\bww_i'$ for $i\in\emptyset$ as the set of all $\bww$ and $\bww'$.)
The assumption of no effects  (Hypothesis \ref{kthorder} with $k=0$) is equivalent to Hypothesis \ref{noeffects}, and the assumption of no  first  and higher order peer effects   (Hypothesis \ref{kthorder} with $k=1$) is equivalent to Hypothesis \ref{nopeer}, 
and the Hypothesis of no  second and higher order peer effects   (Assumption \ref{kthorder} with $k=2$) is equivalent to Hypothesis \ref{nohigher}.

We can also test the hypothesis that there are no direct effects of the own treatment, while allowing for indirect effects from neighbors.
\begin{hypothesis}\label{nohigher}{\sc (No Direct Effects)}
 $Y_i(\bww)=Y_i(\bww')$ for 
all  $i$, and for all pairs of assignment vectors $\bww,\bww'\in\mmw$ such that
$\bww_j=\bww'_j$ for all units $j$ such that $d(i,j)=1$. 
\end{hypothesis}
The most  interesting  version of this null hypothesis might be to test whether the direct effect of the treatment is zero for individuals whose neighbors are all in the control group. This would imply that there could only be a direct effect of the treatment for individuals with at least some treated neighbors. This may be natural in cases where the treatment is some service that requires interacting with other individuals who have the service.

\subsection{Null Hypotheses on Sparsification and Competing Networks}

In the second class of null hypotheses we start with two networks, corresponding to adjacency matrices $\ba_1$ and $\ba_2$. In some cases of interest these may be nested networks, with
$G_{1,ij}\leq G_{2,ij}$ so that $\ba_1$ is a sparsified version of $\ba_2$. 
Suppose we ask individuals whom they regularly interact with, as well whom they have ever interacted with. The first network would define edges using the first question, and the second network would use the second question. 
For example, researchers have used data on emails between employees at Enron to define a network in terms of a threshold for email volume (Goldenberg,  Zheng,
 Fienberg and  Airoldi, 2009).

Alternatively the two networks could correspond to distinct measures of interactions without necessarily being nested, so that for some pairs $(i,j)$, we have $ G_{1,ij}> G_{2,ij}$ whereas for other pairs $(i',j')$, we have
$G_{1,i'j'}<G_{2,i'j'}$. For example, one network definition may be based on email interactions, where another network definition is based on instant messaging interactions, or face-to-face interactions.

We consider the null hypothesis that there is no effect on unit $i$ of the exposure of unit $j$ if $i$ and $j$ are neighbors in the second network $\ba_2$, while allowing for effects on the outcome for unit $i$  of exposure for units $j$ to whom unit $i$ is a neighbor in the first network $\ba_1$.

\begin{hypothesis}\label{sparsehyp}
 $Y_i(\bww)=Y_i(\bww')$ for all $i$, and for  
 all pairs of assignment vectors $\bww,\bww'\in\mmw$ such that
$w_j=w'_j$ for all units $j$ such that $G_{1,ij}=1$. 
\end{hypothesis}

\subsection{Null Hypotheses on Peer Effect Heterogeneity}

Many models of peer effects assume not only that only direct neighbors can influence an individual's outcomes, but also that for any individual it is only the number of treated neighbors that matter, not which of their neighbors got treated. In other words, if we take an individual $i$ with  two neighbors, $j$ and $j'$, the outcome for individual $i$ given assignment $\bww$ with $(\bww_j=0,\bww_{j'}=1)$ is the same as the outcome given  assignment
$\bww'$ with $(\bww'_j=1,\bww'_{j'}=0)$.
Such hypotheses are maintained in many structural models of peer effects, for example the linear-in-means models considered in Manski (1993, 2013).

 Formally:
\begin{hypothesis}\label{noheterohyp}{\sc (No Peer Effect Heterogeneity)}
 $Y_i(\bww)=Y_i(\bww')$ for 
all  $i$, and for all pairs of assignment vectors $\bww,\bww'\in\mmw$ such that
$\sum_{j=1}^N \bww_j\cdot G_{ij}=\sum_{j=1}^N\bww'_j\cdot G_{ij}$. 
\end{hypothesis}
An interesting alternative hypothesis could be that in terms of their effect on outcomes for individual $i$, high-degree neighbors of $i$  are more or less influential than low-degree neighbors of $i$.
This hypothesis implies no second and higher order peer effects, but it is stronger than that. It restricts the range of first order peer effects that is allowed.

A related hypothesis implies that all that matters is that at least one neighbor is exposed to the treatment, and that treating additional neighbors does not affect an individual's outcome. 
\begin{hypothesis}\label{thresholdhyp}{\sc (Threshold Peer Effects)}
 $Y_i(\bww)=Y_i(\bww')$ for 
all  $i$, and for all pairs of assignment vectors $\bww,\bww'\in\mmw$ such that
${\bf 1}\left\{\sum_{j=1}^N \bww_j\cdot G_{ij}>0\right\}={\bf 1}\left\{\sum_{j=1}^N\bww'_j\cdot G_{ij}>0\right\}$. 
\end{hypothesis}
Here an interesting alternative hypothesis could be the number of treated neighbors matters.

\section{Randomization-based Exact P-values with Non-sharp Null Hypotheses: Artificial Experiments}
\label{artificial}

This section contains the main conceptual contribution of the paper.
We describe at an abstract level our approach to the problem of non-sharp null hypotheses. This solution is based on analyzing an artificial experiment that differs from the experiment actually conducted. The artificial experiment is chosen to satisfy two conditions. First, it is chosen so that the original null hypothesis, which was not sharp for the original experiment, is sharp for the artificial experiment, and second, it is chosen so that the randomization-based analysis of the artificial experiment is validated by the design of the original experiment.

We start with an experiment $\cale$, consisting of a set of values  $\mmw$ for the assignment $\bw$, a population  $\mmp$ with $N$ units, and an assignment mechanism $p:\mmw\mapsto [0,1]$. Although in our applications the set $\mmw$ has the structure $\mmw=\{0,1\}^N$, this need not be the case in general. In addition we have a null hypothesis $H_0$ that places restrictions on the function $\by:\mmw\mapsto \mmy^N$.
Instead of testing $H_0$ with the data from this experiment using the randomization distribution implied by $p(\cdot)$, we will analyze a different, artificial, experiment, for which the randomization-based analysis is  validated by the design of  the original experiment. 
Let the artificial experiment be denoted by $\cale'$.
The difference between the artificial experiment and the original experiment has three components. Only one is a choice of the researcher; the remaining two follow from the combination of that choice, the original experiment, and the null hypothesis of interest.

In general test statistics are functions $T:\mmy^N \times \mmw \times \mmx^N\times \mmg\mapsto \mmr$, which are evaluated at $(\by(\bw),\bw,\bx,\ba)$. The first step is to restrict the population whose outcomes the test statistic is allowed to vary with. We denote this subpopulation by $\mmp_F$, and refer to the individuals in this subpopulation as the {focal} units, with $F_i$ an indicator that is equal to one for focal units and zero otherwise. In the special case where the null hypothesis is that of no spillovers at all, the focal subpopulation corresponds to the subpopulation of fixed units in Aronow (2012), who refers to its complement  as the variant units. However, because in our approach the artificial experiment may also need to hold fixed the treatment assignment for some units outside the subpopulation of   what Aronow calls the fixed subpopulation, we use a different terminology. At this point the choice of focal subpopulation is arbitrary. Its choice does not affect the validity of the resulting p-values, but as we shall discuss below, it has  a major impact on  the power of the test.
Let $N_F$ be the cardinality of the set $\mmp_F$,  let $\by_F(\bww)$ denote the $N_F$-vector of potential outcomes for the focal units for any  treatment $\bww$, and let $\by^\obs_F=\by_F(\bw)$ be the vector of realized outcomes for these units given the actual assignment $\bw$. 
The selection of this subpopulation can depend generally on the fixed characteristics of the population  $\bx$, and the network $\ba$. It cannot depend on the assignment $\bw$ either directly, or indirectly through dependence on the realized outcome $\by^\obs$. We now consider test statistics $T:\mmy^{N_F}\times\mmw\times\mmx^N\times \mmg\mapsto \mmr$, evaluated at $(\by_F(\bw), \bw, \bx,\ba)$.


Given the focal subpopulation $\mmp_F$ and the null hypothesis $H_0$, define the set of subsets of $\mmw$,
\[ \mms=\cup_{w\in\mmw} \Bigl\{\cap_{i\in\mmp_F} \mmv(i,\bww,H_0)\Bigr\}.\]
This set plays a key role in our approach. An important property  is that it is a partition of $\mmw$.
\begin{proposition}{\sc (Partition of the Assignment Space)}
$\mms$ is a partition of $\mmw$.
\end{proposition}
{\bf Proof:} Because $\bww\in\cap_{i\in\mmp_F}\mmv(i,\bww,H_0)$, it immediately follows that $\cup_{\mmv\in\mms}\mmv=\mmw$. Thus the remaining property to be established is that either $\left(\cap_{i\in\mmp_F}\mmv(i,\bww,H_0)\right)\cap \left(\cap_{i\in\mmp_F}\mmv(i,\bww',H_0)\right)=\emptyset$ or $\cap_{i\in\mmp_F}\mmv(i,\bww,H_0)= \cap_{i\in\mmp_F}\mmv(i,\bww',H_0)$. 
If $\left(\cap_{i\in\mmp_F}\mmv(i,\bww,H_0)\right)\cap \left(\cap_{i\in\mmp_F}\mmv(i,\bww',H_0)\right)$ is not equal to the empty set, there must be a 
$\bww''\in \mmv(i,\bww,H_0)$ and $\bww''\in\mmv(i,\bww',H_0)$. Then 
\begin{equation}\by_F(\bww'')=\by_F(\bww')=\by_F(\bww).\label{een}\end{equation}
Hence if there is another element $\bww'''\in \mmv(i,\bww',H_0)$, it must be the case that
\[ \by_F(\bww''')=\by_F(\bww').\]
By (\ref{een}) this is equal to $\by_F(\bww'')$, and also be (\ref{een}) this is equal to $\by_F(\bww)$. Hence
 it must be the case that
\[ \by_F(\bww''')=\by_F(\bww'')=
\by_F(\bww')=\by_F(\bww),\]
and $\bww'''\in \mmv(i,\bww,H_0)$. Therefore $\cap_{i\in\mmp_F}\mmv(i,\bww,H_0)= \cap_{i\in\mmp_F}\mmv(i,\bww',H_0)$, which finishes the proof.
$\square$.

The third component of the artificial experiment consists of a new assignment mechanism $p':\mmw\mapsto[0,1]$.
To define this third component we decompose the original experiment into a stratified experiment. Given the partition $\mms$, define
the stratum indicator $S:\mmw\mapsto \{1,\ldots,J\}$, so that the stratum is $S(\bww)=j$ if $\bww\in\mmw_j$.
 Now we can think of the original experiment $\cale$ as a stratified experiment where we first draw the stratum $S$, with $\pr(S=j)=\pr(\bw\in\mmw_j)$, followed by the second stage where we draw $\bw$  conditional on $S$, with 
\[p'(\bww)=\pr(\bw=\bww|S=j)=
\left\{\begin{array}{ll}\frac{p(\bww)}{\sum_{\bww'\in\mmw_j} p(\bww')},\hskip0.5cm & {\rm  if }\ \pr(S=j)>0,\bww\in\mmw_S,\\
 0 & {\rm otherwise.}\end{array}\right.\]


Now we propose to analyze the artificial experiment $\cale'=(\mmw_S,\mmp_F,p'(\cdot))$. The set of restrictions on the values the function $\by:\mmw\mapsto \mmy^N$ that corresponds to the original null hypothesis translates into a set of restrictions on the values of the function $\by_F:\mmw_S\mapsto\mmy^{N_F}$ which gives us the implicit null hypothesis for the new experiment. By contstruction, the set of assignments $\mmw$ and the focal population $\mmp_F$ are chosen so that the null hypothesis is sharp for this artificial experiment. Formally, for any pair $(\bww,\by_F(\bww))$ with $\bww\in\mmw_S$, we can infer the values of $\by_F(\bww')$ for any other value $\bww'\in \mmw_S$. We discuss some examples of this in the next section.
 We then choose a statistic $T: \mmy^{N_F} \times \mmw \times \mmx^N\times\mmg\mapsto \mmr$ that depends only on the  outcomes for the individuals in the focal population, $\by^\obs_F=\by_F(\bw)$. We calculate p-values for this statistic by comparing the realized value of the statistic, $T^\obs=T(\by_F^\obs, \bw,\bx,\ba)$, to the  randomization distribution  for $T(\by_F(\bw),\bw,\bx,\ba)$ induced by the modified assignment distribution $p'(\cdot)$. 

A key insight is that a randomization-based analysis of the artificial experiment $\cale'$ is validated by the design of the original experiment $\cale$. Let us consider the two modifications--changing the population and using a conditional assignment mechanism--in turn and justify this claim. Choosing a subpopulation of units based on fixed attributes or pretreatment variables  such that the test statistic varies only with outcomes for these units does cannot invalidate the p-value because it is valid for any statistic.
Second, consider the change in the assignment mechanism. We can think of the original assignment mechanism, corresponding to the distribution $p(\cdot)$, as a two-stage procedure: first we choose $S$, and then the actual assignment is determined either by drawing according to $p'(\cdot)$ where  $p'(\bww)={\rm pr}(\bw=\bww|\bw\in\mmw_S)$.
Thus the artificial experiment conditions on the value of $S$ and only exploits the second stage randomization. In general this may discard information, but it does not affect the validity.

\section{Exact P-values for the Null Hypothesis of No Spillovers}\label{nospillover}

Here we discuss how exact p-values can be calculated for the hypotheses introduced in Section \ref{hypotheses}, given randomized assignment of the treatments. To simplify the discussion we focus in this section initially on  a completely randomized experiment, where $M$ units out of $N$ are randomly selected to receive the treatment (see Imbens and Rubin, 2015 for a general discussion). In Section \ref{cluster} we discuss extensions to clustered randomized experiments.
\begin{assumption}\label{randomassignment}{\sc (Random Assignment)}
\[ {\rm pr}(\bw=\bww)=1\biggl/\left(\begin{array}{c}N \\ M\end{array}\right),\]
for all $\bww \in \{0,1\}^N$ such that $\sum_{i=1}^N \bww_i=M$.
\end{assumption}

To set the stage, let us first consider the case where we test the null hypothesis of no treatment effects whatsoever. In that case for each individual  $\mmv(i,\bww,H_0)=\mmw$, we can take the subpopulation of focal units to be the entire population, $\mmp_F=\mmp$,  and the partitioning is $\mms=\{\mmw\}$. Then the assignment mechanism is the same under the artificial experiment as it is under the original experiment, $p'(\cdot)=p(\cdot)$, and thus the artificial experiment is identical to the original experiment.

\subsection{Exact P-values for the Null Hypothesis of No Spillovers when the Network consists of Dyads}

To develop some intuition for the problem we first look at the case where the network has a simple structure.
Suppose the population consists of $N$ units paired into $N/2$ dyads. For individual $i$ let $\ell(i)\in\{1,\ldots,N\}$ be the index of the neighbor of individual $i$.  We are interested in testing the hypothesis that there are no spillover effects (Hypothesis \ref{nopeer}),  allowing for the possibility of direct effects of the own treatment on an individual's own outcome. 

\subsubsection{The Artificial Experiment}

To create the artificial experiment $\cale'$ we first select the focal subpopulation. We do this by selecting one member from each pair, and designate that individual in the pair as the focal individual. 
This selection can be random, or based on pretreatment variables, but not on outcome or assignment data.
Let $F_i=1$ if an individual is a focal individual  and $F_i=0$ for non-focal, or auxiliary individuals. Selecting one focal unit from each pair is not required for our approach, but it makes intuitive  sense. If both members of a pair are focal units, then the level sets imply that we cannot vary the treatments for any member of the pair in the artificial experiment. If neither member of the pair is focal we do not use the outcomes for the two units. In both cases the pair is essentially dropped from the analysis, so only if there is a single focal unit in each pair does the pair enter in the analysis.

In the second step, we define the restricted set of assignments $\mmw_S$.
Let $\bw$ be the full assignment vector.
For individual $i$, $\mmv(i,\bww,H_0)=\{\bww'\in\mmw|\bww_i=\bww'_i\}$.
Hence
\[ \mmw_S=\cap_{i\in\mmp_F} \mmv(i,\bw,H_0)=
\{\bww\in\mmw|\bww_i=W_i\ {\rm for\ all}\ i\in\mmp_F\},\]
allowing only the treatments for the non--focal, or auxiliary units, to vary.
Let $M_F=\sum_{i:F_i=1} W_i$ be the number of treated focal individuals, and $M-M_F$ the number of treated auxiliary individuals. Then, because there are $N/2$ auxiliary individuals, the distribution of assignments $p'(\cdot)$ in the artificial experiment satisfies
\[ p'(\bww)=
{\rm pr}(\bw={\bww}|S)=1\biggl/\left(\begin{array}{c}N/2 \\ M-M_F\end{array}\right),\]
 for $\bww\in\mmw_S$, and zero otherwise.

Given the experiment we consider test statistics $T:\mmy^{N_F} \times \mmw \times\mmx\times\mmg\mapsto \mmr$. 
For any statistic in this class we can infer its distribution under the null hypothesis. We would like to choose the statistic whose distribution is sensitive to  interesting departures from the null hypothesis. 
We consider two statistics, motivated by  parametric models that allow for spillover effects.

\subsubsection{Test Statistics}

Consider a model for the potential outcomes that does not impose the null hypothesis of no spillovers. In that case, with a single neighbor for each individual, the potential outcome for individual $i$ can be written as a function of the own treatment $w_i$ and the neighbor's treatment $\bww_{\ell(i)}$, or, $Y_i(\bww)=Y_i(\bww_i,\bww_{\ell(i)})$. A natural starting point is to assume that both direct (own) and indirect (neighbor's) treatment effects are constant and additive:
\begin{equation}\label{addmodel}
 Y_i(\bww_i,w_{\ell(i)})=\alpha+\tau_{\direct}\cdot \bww_i+
\tau_\spill\cdot \bww_{\ell(i)}+\varepsilon_i.\end{equation}
Given this parametric model the null hypothesis of no spillovers  corresponds to $\tau_\spill=0$.
To find a statistic with good power properties for testing our nonparametric null hypothesis of no-spillovers, we can look at the Lagrange multiplier or score  test statistic for the  null hypothesis $\tau_\spill=0$ in this parametric model, assuming homoskedasticity, normality and independence for the $\varepsilon_i$. The validity of our proposed testing procedure does not rely on these parametric and distributional assumptions, but if they hold, the fact that in that case the test corresponds to a Lagrange multiplier test  would endow the procedure with large sample efficiency.

In this parametric model the likelihood function for the focal units
is 
\[{\cal L}(\sigma^2,\alpha,\tau_\direct,\tau_\spill)
=\prod_{i:F_i=1} \frac{1}{\sqrt{2\pi\sigma^2}}\exp\left(-\frac{1}{2\sigma^2}
\left(Y^\obs_i-\alpha-W_i\cdot\tau_\direct-W_{\ell(i)}\cdot \tau_\spill\right)^2\right),\]
where $\sigma^2$ is the variance of $\varepsilon_i$.
The sum of the scores, that is, the  sum over the focal units of the derivative of the logarithm of the density under this model with respect to $\tau_\spill$, evaluated at $\tau_\spill=0$, is equal to
\[ {\cal S}=\frac{1}{\sigma^2}\sum_{i=1}^N W_{\ell(i)} \cdot \Bigl(
Y^\obs_i- \alpha-\tau_{\direct}\cdot W_i\Bigr).\]
The statistic we focus on is this sum with $\alpha$ and $\tau_\direct$ replaced by estimates based on the outcomes for only the focal units. These estimates are
\[ \hat\alpha=\oy^\obs_{F,0},\hskip1cm 
\hat\tau_\direct=\oy^\obs_{F,1}-\oy^\obs_{F,0},\]
where, for $w=0,1$, $\oy^\obs_{F,w}$ is the average outcome for focal units with $W_i=w$ and $N_{F,w}$ is the number of focal units with $W_i=w$.
This leads to the statistic, after normalizing by the number of focal units,
\begin{equation}\label{dyad1} T_{\text{score}}^{\text{dyad}}=\frac{1}{N_F}\sum_{i:F_i=1}  \Bigl(
Y^\obs_i-\oy^\obs_{F,0}-W_i\cdot\left(\oy^\obs_{F,1}-\oy^\obs_{F,0}\right)
\Bigr)
\cdot W_{\ell(i)}
.\end{equation}

This statistic is interpreted as the correlation between the neighbors' treatment status and the focal unit's outcome, adjusted for the average value of the outcome
for focal units with the same treatment status. 

Although such a model appears substantively less plausible, it is also interesting to consider the model in (\ref{addmodel}) without a direct effect:
\begin{equation}\label{addmodel2} Y_i(w_i,w_{\ell(i)})=\alpha+
\tau_\spill\cdot w_{\ell(i)}+\varepsilon_i.\end{equation}
Then the Lagrange multiplier approach leads to the statistic
\begin{equation}\label{dyad2} T_{\text{elc}}^{\text{dyad}}=\frac{1}{N_F}\sum_{i:F_i=1} W_{\ell(i)}\cdot \Bigl(
Y^\obs_i-\oy^\obs_F
\Bigr)=\frac{N_{F,(1)}}{N_{F,(0)}}\cdot\Bigl(\oy^\obs_{F,(1)}-\oy^\obs_{F,(0)}\Bigr),\end{equation}
where
for $w=0,1$, $\oy^\obs_{F,(w)}$ is the average outcome for focal units with neighbors whose treatment status is $W_{\ell(i)}=w$ and $N_{F,(w)}$ is the number of focal individuals whose neighbor has treatment status $w$.  
Hence the statistic essentially compares average outcomes for focal units with treated neighbors and focal units with control neighbors.  We refer to this statistic as an 
edge-level-contrast statistic for reasons that will become clear below when we generalize the network structure.

The first statistic, $T^{\text{dyad}}_{\text{score}}$, yields a more powerful test when there are direct effects of the treatment, because it adjusts for the estimated direct effects of treatment.
Failing to do so introduces additional noise in the distribution of the test statistic.

\subsection{Artificial Experiments for the Null Hypothesis of No Spillovers for General Networks}

In this section we consider the more general problem of testing for spillover effects in an unrestricted network setting. We maintain the assumption that the randomization is at the unit level, with $M$ randomly selected units out of the population of $N$ units exposed to the intervention.
As before we choose a subpopulation of focal individuals  whose outcomes we use, with the complement of this subpopulation the set of auxiliary individuals. This selection may be random or depend on pretreatment variables.
The restricted set of assignments  fixes the assignments for the focal individuals:
$ \mmw_S=\{\bww\in\mmw|\bww_F=\bw^\obs_F\},$ allowing only the treatments for the non-focal or auxiliary units to vary.
There are two substantive differences 
 with the setting where the network consists of pairs. The choice of the statistic is more complicated, and so is the choice of the focal subpopulation.

\subsection{Test Statistics}

We consider three test statistics.  The first is a modification of a test statistic previously proposed by Bond et al. (2012); the second is optimal for a particular data-generating process; and the third is a modification of a statistic proposed by Aronow (2012).

\subsubsection{The Edge-Level Contrast Statistic}\label{bond}

The first  statistic we consider is a modification of an edge-level statistic used by Bond et al. (2012). Bond et al. test for the presence of spillovers using the randomization distribution based on the null hypothesis of no effects of the treatment whatsoever. The statistic they use is equal to the difference between the average, over all edges where 
the alter is exposed to the treatment, of the ego's outcome and the average, over all edges where the alter is not exposed to the treatment, of the ego's outcome: 
\[ T_{\rm B}(\bw,\by^\obs,\ba)=\frac{\sum_{i,j\neq i} G_{ij}\cdot W_j\cdot Y_i^\obs}
{\sum_{i,j\neq i} G_{ij}\cdot W_j}
-\frac{\sum_{i,j\neq i} G_{ij}\cdot (1-W_j)\cdot Y_i^\obs}
{\sum_{i,j\neq i} G_{ij}\cdot (1-W_j)}.\]
We cannot infer the randomization distribution of this statistic if we only impose the null hypothesis of no spillovers  but allow for direct effects of the treatment (which is the null hypothesis of interest). Bond et al. report p-values based on the additional assumption that there are no own effects of the treatment.
Without this additional assumption the p-values reported based on this statistic are therefore not generally valid. In Appendix A we provide analytical calculations  that show that the size distortions for this statistic can be substantial in the presence of direct effects of the treatment, as high as 0.2 for a nominal 0.05 level test in simple cases.

However,  we can  modify the Bond et al. statistic,  averaging only over the subset of edges where the ego is in the focal subpopulation and the alter is in the auxiliary subpopulation (in the current setting where we test the null of spillovers this subpopulation is equal to the complement of the focal subpopulation):
\begin{equation}\label{mbf} T_{\elc}(\bw,\by^\obs_F,\ba)\end{equation}
\[\hskip1cm=\frac{\sum_{i,j\neq i} F_i\cdot G_{ij}\cdot (1-F_j)\cdot W_j\cdot Y^\obs_i}
{\sum_{i,j\neq i} F_i\cdot G_{ij}\cdot (1-F_j)\cdot W_j}
-\frac{\sum_{i,j\neq i} F_i\cdot G_{ij}\cdot (1-F_j)\cdot (1-W_j)\cdot Y^\obs_i}
{\sum_{i,j\neq i} F_i\cdot G_{ij}\cdot (1-F_j)\cdot (1-W_j)}.\]
We refer to this as the edge-level-contrast statistic.
In the case where the network consists of dyads, it reduces to our second test statistic for the case of dyads, $T^{\rm dyad}_{\text{elc}}$ in (\ref{dyad2}).

\subsubsection{A Score Test Statistic}
\label{score}

We motivate the second test statistic in a more systematic way with a structural model for treatment effects. 
Suppose we use a simple linear model, a simplified version of the  linear-in-means model of the type discussed in Manski (1993, 2013) with only exogenous peer effects:
\begin{equation}\label{lim} Y^\obs_i=\alpha_0+\tau_\direct\cdot W_i+\tau_{\rm exo}\cdot
\sum_{j=1}^N W_j\cdot \overline{G}_{ij}
+\varepsilon_i,\end{equation}
where $\overline{G}_{ij}=G_{ij}/\sum_{j'=1}^N G_{ij'}$ is a normalized indicator for links.
(If $\sum_{j'=1}^N G_{ij'}=0$, then $\overline{G}_{ij}=0$.) 
Hence $\sum_{j=1}^N W_j\cdot \overline{G}_{ij}$ is the fraction of treated friends.

Testing for spillovers in the context of this model corresponds to testing the parametric null hypothesis that the exogenous peer effects parameter $\tau_{\rm exo}$ is equal to zero.
A natural way to derive a powerful test statistic for $\tau_{\rm exo}=0$ in a parametric  model, and the basis of Lagrange multiplier tests, is to derive the average score for $\tau_{\rm exo}$, evaluated at $\tau_{\rm exo}=0$ and estimates for the nuisance parameters ($\alpha_0$ and $\tau_\direct$ in this case).
Under the model in (\ref{lim}) 
  the score statistic is proportional to the covariance between the residual under the null and 
the fraction of neighbors who are treated, $\sum_{j=1}^N \overline{G}_{ij}\cdot W_j$, leading to
\begin{equation}\label{corr} T_{\rm score}=
{\rm Cov}\left(\left.
Y^\obs_i-\hat\alpha-\hat\tau_\direct\cdot W_i,
\sum_{j=1}^N W_j\cdot \overline{G}_{ij}
\right| \sum_{j=1}^N G_{ij}>0,F_i=1
\right)
.\end{equation}
\noindent{\sc Remark 1}
If the network consists of dyads, with one unit in each dyad designated focal and the other auxiliary, then this statistic is identical to the statistic $T^{\rm dyad}_{\text{score}}$ in (\ref{dyad1}).  
As in the case of dyads, this test statistic reduces variance in the test statistic by normalizing outcomes by the estimated direct effect of the treatment, at least when  direct effects of the treatment ar present. $\square$

\noindent{\sc Remark 2}
Note that our approach to deriving the test statistic can be applied to alternative structural models with different functional forms for outcomes, the nature of spillovers, etc., and
as above, the test statistic is valid irrespective of the validity of the structural model.  The 
 power of the test, however, will depend on the quality of the model.
$\square$

\noindent{\sc Remark 3}
It is also interesting to note that the same score statistic applies to a different model. Suppose we start with a different version of the linear-in-means model of the type discussed in Manski (1993, 2013):
\begin{equation}\label{lim} Y^\obs_i=\alpha_0+\tau_\direct\cdot W_i+\tau_{\rm endog}\cdot\overline{Y}^\obs_{(i)}+\varepsilon_i,\end{equation}
where $\overline{Y}^\obs_{(i)}$ is the average outcome for $i$'s neighbors. In this model the spillovers arise from the direct effect of one's own treatment  on one's own outcome (if $\tau_\direct\neq 0$), combined with what Manski calls endogenous effects of the neighbors' outcome on the own outcomes ($\tau_{\rm endog}$). This implies that treatment exposure for non-neighbors can affect one's outcome if the non-neighbor are connected through other individuals, with the magnitude of the spillover effects depending on the distance between the individuals in the network. Although this endogenous peer effects model implies that spillover effects propagate throughout one's network, the score statistic for this model is identical to that in (\ref{lim}), because close to the null of no spillover effects the effects are dominated by those of direct neighbors.
Details for this calculation are presented in Appendix B. $\square$

\subsubsection{The Has-Treated-Neighbor Test Statistic}
As the third test statistic, we consider a variation on a statistic based on distance to the nearest treated unit.
Aronow (2012) proposes a test statistic for spatial or network interference that is the correlation between outcome for focal units and the distance to the nearest treated auxiliary unit.
If distance is defined in terms of hops between two units in a network and there are many treated units, then much of the variation in this measure will be between having a treated unit in one or two hops.
So we analyze a related statistic the uses, instead of the distance to the nearest treated unit, an indicator for whether any of a unit's non-focal neighbors are treated. This statistic is the correlation between this indicator and the outcome, both for focal units:
\[ T_{\rm htn}=\frac{1}{S_{Y^\obs_F}\cdot S_{\rm TA}}\frac{1}{N_F}\sum_{i\in \mmp_F}
\left(Y^\obs_i-\overline{Y}_F^\obs\right)
\cdot {\bf 1}_{\sum_{j} G_{ij}\cdot W_j\cdot (1-F_j)>0},\]
where $S_{Y^\obs_F}$ and $S_{\rm TA}$ are the sample standard deviation of the outcome for focal units and the standard deviation for the indicator, for focal units, of having at least one treated auxiliary neighbor.  Like the edge-level contrast statistic, this statistic does not adjust for estimated direct effects of the treatment.

\subsection{Choosing the Focal Subpopulation for the Null Hypothesis of No Spillovers}

\label{focal}

A key feature of our approach is that the researcher needs to choose a focal subpopulation. This choice, in combination with the null hypothesis, determines the  randomization distribution in the artificial experiment. Although the p-values are valid irrespective of the choice of focal subpopulation, this choice may affect the power of the testing procedure substantially.

Here we discuss some algorithms for choosing the subpopulation of focal units, where the goal is to maximize the power of the test. In general the power will depend on a number of features of the problem. First, it will depend the alternative hypothesis, for example whether the spillover effects are linear in the number or the proportion of treated neighbors. Second, the power will depend on the choice of statistic. The power will also depend on the network structure. Finding the focal subpopulation that optimizes power for particular choice of alternative and a particular test statistic is a difficult problem. Here we discuss some issues and suggest general solutions that may have good power in a wide range of settings.

In the case of testing the null of no spillovers, there are three general principles that apply irrespective of the specific alternative hypothesis and test statistic. First, because the artificial experiment considers only change in the treatment for auxiliary individuals, it is important that there are a substantial number of auxiliary individuals. Second, because the statistic depends only on outcomes for focal units, it is important that there is a substantial number of focal units. Third, because the alternative hypothesis involves spillovers from treated alters to focal egos, and because only changes in the treatment for auxiliary individuals are considered, it is important that there are many edges between focal and auxiliary individuals. These principles were helpful in the dyad case, where they suggested selecting a single focal individual in each pair.
Some settings may also have additional constraints that guide the selection of focal units. For example, we might only observe the outcome for a small fraction of the units even though the treatment is observed for all units (e.g., Bond et al. (2012) only observe voting status for about 10\% of their population).

\subsubsection{Random Selection}

As a baseline method we randomly choose 50\% of the population to be focal, with the remainder auxiliary, without regard to the network structure.

\subsubsection{Selection Based on $\varepsilon$-Nets}
In the second approach to focal unit selection, we aim to select a large set of focal units that are not adjacent to each other. In particular, we use a method for finding an $\varepsilon$-net (see, e.g., Gupta, Krauthgamer and Lee, 2003), or a set of points that is both an $\varepsilon$-packing and an $\varepsilon$-covering, with $\varepsilon = 2$.\footnote{A 2-net is also called an \emph{independent set} and the greedy algorithm we give here constructs a maximal independent set. We describe this in terms of $\epsilon$-nets because larger values of $\epsilon$ might be used when testing other hypotheses about spillovers.}
To define an $\varepsilon$-net on a graph, we let $B_\varepsilon(i) = \{j : d(i, j) \leq \varepsilon$ and $j \in \mmp\}$ be the set of all vertices within $\varepsilon$ hops of vertex $i$.

\begin{definition}\label{def_epsilonnet}{\sc ($\varepsilon$-net in a graph)} An $\varepsilon$-net is a set of vertices $\mms \subseteq \mmp$ such that: (a) the vertices are mutually at distance at least $\varepsilon$ from each other, $d(i, j) \geq \varepsilon$ for all $i, j \in \mms$; and (b) the union of all of their $\varepsilon$-balls covers all vertices,
$\mmp \subseteq \cup_{i \in \mms} B_\varepsilon(s)$.
\end{definition}
Ugander, Karrer, Backstrom, and Kleinberg (2013) describe a greedy method for finding an 3-net, which can be generalized to find a $\varepsilon$-net for other values of $\varepsilon$. To find a 2-net, we do the following.
Starting with an empty set of focal units and an empty set of auxiliary units we randomly select
 a seed for the $\varepsilon$-net. Given the new seed we assign it to the focal subpopulation, and we assign all of its neighbors to the auxiliary subpopulation. If at that point all individuals are assigned to either the focal or the auxiliary subpopulation we stop. If not, we randomly draw another seed to be assigned to the focal subpopulation and assign all its neighbors to the auxiliary subpopulation. We continue randomly selecting new seeds until all individuals are assigned to either the focal or auxiliary subpopulation. This greedy algorithm leads to a set of focal units that are not neighbors.

\subsubsection{Maximizing the Number of Edge Comparisons}
\label{imbensnet}


In the third approach we choose the focal subpopulation by attempting to maximize the number of focal--auxiliary edges, 
\[ N(\bff,\ba)=\sum_{i,j} F_i\cdot  G_{ij}\cdot (1-F_j),\]
leading to\[ \bff^*=\arg\max_\bff N(\bff,\ba).\]
  This approach ignores the fact that the average over the edges may involve multiple edges with the same ego. This would not change the optimality if the number of focal-auxiliary edges were the same for all focal individuals, but if there is substantial variation in the number of such edges one might do better taking that into account.

Solving this problem exactly is computationally demanding, so we approximate it by using a greedy algorithm.
We start by assigning all units to the auxiliary subpopulation, so that there are no focal-auxiliary edges.
We then calculate for each non-focal unit the number of focal-auxiliary edges that would get added if unit $i$ gets moved to the focal subpopulation,
$\Delta_{N,i}$. Next, add
the individual to the focal subpopulation who bring the biggest gain. This process continues until no additional focal unit would increase the number of focal-auxiliary edges.

Suppose we have an initial focal subpopulation $\bff$. For  auxiliary individual $i$ consider adding them to the focal subpopulation. That would change $N(\bff,\ba)$ by the number of the auxiliary neighbors of $i$  minus  the number of  focal neighbors of $i$: 
\[ \Delta_{N,i}=K_{A,i}- K_{F,i}.\]
This puts a premium on selecting focal units with a larger number of edges. Because we consider settings where it is the fraction of neighbors that are treated that matters for the spillover effects, rather than the total number, we modify this criterion by dividing it by the number of neighbors, and selection as an additional focal unit the one with the highest value for 
\[ \delta_{N,i}=\frac{K_{A,i}- K_{F,i}}{K_i}.\]
In regular graphs (i.e., where all units have the same number of neighbors) this change does not matter, but it does in settings with where the degree distribution has a positive variance. Thus, we sequentially add to the set of focal units the unit $i$, among those currently not in the focal subpopulation, who has the highest value for $\delta_{N,i}$, until there is no auxiliary unit with a positive value for $\delta_{N,i}$.

In settings where the network consists of dyads, both the $\varepsilon$-net approach and maximizing the number of edge comparisons leads to the same result: in each dyad one randomly selected vertex will be the focal unit and the other vertex in the dyad will be the auxiliary unit. In that case the random selection of focal units without regard to network structure will be substantially less powerful by allowing for the possibility that both individuals in a dyad are focal or that both are auxiliary.

There are more general connections between this method and the 2-net method. With the modified, fractional criterion $\delta_{N,i}$, this method first selects a 2-net as the focal units and then continues to add focal units. That is, this method allows using a larger set of focal units than would be selected by finding a 2-net.

\subsection{Exact P-values for Spillovers with Clustered Random Assignment}
\label{cluster}

Now suppose the randomization is more complex than the one considered in the previous section, where we randomly selected $M$ units out of the population of $N$ to receive the treatment. Of particular interest is the generalization with clustered randomization. In this case the population is first partitioned into $K$ clusters, $\mmp_1,\ldots,\mmp_K$, with $\mmp_k\subset\mmp$, $\mmp_k\cap\mmp_l=\emptyset$ if $k\neq l$, and $\cup_{k=1}^K\mmp_k=\mmp$. This partitioning may depend on the network structure. 
In fact, in graph cluster randomization, the partitioning is often chosen so as to heuristically maximize the fraction of edges within that are within clusters, subject to other constraints (e.g., cluster size), or other related quantities, such as modularity (Newman, 2006). See Eckles, Karrer, and Ugander (2014) and Ugander, Karrer, Backstrom, and Kleinberg (2013).
Let $C_i\in\mmc=\{1,\ldots,K\}$ indicate the cluster that individual $i$ belongs to.
In the next step, $M$ of the $K$ clusters are assigned to the treatment group, implying all units in those $M$ clusters will be exposed to the treatment, and the remaining units will be assigned to the control group. More generally, we may consider an unrestricted distribution for the assignment vector ${\bf W}$, specified by the function
$p:\mmw\mapsto [0,1]$ for some set of assignments $\mmw$ that is different from one that assigns equal probability to all assignments with $M$ treated and $M-N$ control units. 

For the original experiment the clustering does not change the fundamental approach. If we are interested in testing a sharp null hypothesis such as the null hypothesis of
no effect of the treatment whatsoever, we can use exactly the same statistics. The only difference is that when we calculate the distribution of the statistic under the null, we now do so under the assignment mechanism defined by the clustering. Because many assignment vectors $\bww$ that are possible under complete randomization are ruled out under cluster randomization, the clustering typically reduces the power of the  tests.
This issue is even more of a concern for testing null hypotheses regarding spillovers.
We again select a focal subpopulation $\mmp_F\subset\mmp$. For each individual calculate the set of assignments that do not change the outcome for that individual under the null hypothesis,  $\mmv(i,\bww,H_0)$. The restricted set of assignments is, as in the general case, the intersection of these  sets over all focal individuals:
\[ \mmw_S=\prod_{i\in\mmp_F} \mmv(i,\bw,H_0).\]
The  distribution of the assignments in the artificial experiment is, as before, the conditional probability given that $\bw\in\mmw_R$:
\[ p'(\bww)=\frac{p(\bww)}{\sum_{\bww'\in\mmw_S} p(\bww')},\]
for $\bww\in\mmw_S$, and zero elsewhere.
The artificial experiment is now characterized by the triple $(\mmw_S,\mmp_F,p'(\cdot))$.

For any statistic
$ T:\mmw_S\times \mmy^{N_F}\times\mmx\times\mmg\mapsto \mmr,$
we can infer its exact distribution under the null hypothesis of no spillovers, using the randomization distribution induced by the clustered randomization. Thus we can use the same statistics as before, e.g., the edge-level-contrast statistic or the score statistic. The change in the distribution of the treatment affects the power of the tests, but does not fundamentally change the approach.

To illustrate what practical issues the clustered randomization raises, consider the edge-level-contrast statistic $T_{\rm elc}$. This statistic is equal to the difference in the average outcome for focal units over all edges between one focal unit and one auxiliary unit, where the auxiliary unit is treated and the average outcome for focal units over all edges where the auxiliary unit is in the control group. Because treatments for units in the same clusters as  focal units do not vary in $\mmw_S$ because of the cluster randomization, the power of the tests will be severely reduced if the clusters are constructed in such a way that there are few between-cluster edges. Although such clustering designs may be effective in estimating total causal effects that include both direct effects and spillover effects, e.g., Eckles, Karrer, and Ugander (2014) and Ugander, Karrer, Backstrom, and Kleinberg (2013), they may be less suited towards distinguishing between the two effects.

\section{Exact P-values for the Null Hypothesis of No Higher Order Peer Effects}\label{higher}

Now consider the case where we are interested in the null hypothesis of no higher order peer effects, Hypothesis \ref{nohigher}. We focus again on the case with complete random assignment, although that is not critical. Define $\bh$ to be the matrix indicating neighbors of neighbors, so that 
\[H_{ij}=\left\{
\begin{array}{ll}
1\hskip1cm & {\rm if}\ i\neq j\wedge G_{ij}=0\wedge\Bigl(
\sum_{k=1}^N G_{ik}\cdot G_{jk}>0\Bigr)\\
0 & {\rm otherwise.}
\end{array}\right.
\]
Again select a focal subpopulation $\mmp_F$. The change in the null hypothesis does not impose restrictions on the choice of the focal subpopulation, although the implications of this choice for the power are different compared to the case where the null hypothesis ruled out the presence of any spillovers. 
The difference with the previous null hypothesis of no spillovers is in the definition of the restricted set of assignments $\mmw_S$.
Given this null hypothesis, for individual $i$, the level set $\mmv(i,\bww,H_0)$ now consists of the set of assignments $\bww'$ such that the assignments are the same for $i$ and for all $i$'s neighbors
\[ \mmv(i,\bww,H_0)=\{\bww'\in\mmw|\bww'_i=\bww_i\wedge \left(\bww'_j=\bww_j{\rm \ for\ all}\ j \ {\rm s.t.}\ G_{ij}=1\right)\}.\]
Then, as before, the restricted set of assignments is the intersection over all focal units of these sets:
\[ \mmw_S=\cap_{i\in\mmp_F} \mmv(i,\bw,H_0).\]
We can conceptualize this set in terms of a partition of the population into three subpopulations. Given the subpopulation of focal units $\mmp_F$, define the set of buffer units $\mmp_B$ who are not focal, but who have one or more neighbors who are focal:
\[ \mmp_B=\left\{i\in\mmp\left| F_i=0\wedge \left(\sum_{j=1}^N G_{ij}\cdot F_j>0\right)\right.\right\},\]
and the set of auxiliary units $\mmp_A$ who are not focal, nor do they have neighbors who are focal:
\[ \mmp_A=\left\{i\in\mmp\left| F_i=0\wedge \left(\sum_{j=1}^N G_{ij}\cdot F_j=0\right)\right.\right\}.\]
Then the restricted set of assignments keeps fixed the assignment for units who are not auxiliary, that is, for focal and buffer units:
\[ \mmw_S=\{\bww\in\mmw|\bww_i=\bw_i{\rm \ if\ } i\in\mmp_F\cup\mmp_B\}.\]

To visualize this consider a very simple example with a population with three units, with the only edge between individuals 1 and 2, corresponding to the following adjacency matrix:
\[ \ba=\left(\begin{array}{ccc}
0 & 1 & 0 \\
1 & 0 & 0\\
0 & 0 & 0\\
\end{array}
\right).\]
Suppose we choose unit 1 to be the focal unit, $\mmp_F=\{1\}$. Then
 the set of neighbors of focal units, or the set of buffer units
 is $\mmp_B=\{2\}$ and  the set of auxiliary units is $\mmp_F=\{3\}$.
Suppose the actual assignment is $\bw=(0,0,0)$.
Then
\[ \mmw_S=\mmw(1,\bw,H_0)=\{(0,0,0),(0,0,1)\},\]
allowing only the assignments for the auxiliary unit to vary.

Now, the experiment we consider is that of randomly assigning $\bw$ within the set $\mmw_S$. Under those assignments we know all the potential outcomes for focal individuals.
The new assignment mechanism is, as before, the conditional assignment probability given the assignments for non-auxiliary units,
$ p'(\bww)={\rm pr}(\bw=\bww|\bw\in\mmw_S),$
and the artificial experiment is
\[ \cale'=(\mmw_S,\mmp_F,p'(\cdot)).\]

\subsection{Test Statistics}

Let us now consider test statistics for this setting.

\subsubsection{An Edge-Level-Contrast Statistic}

A natural approach to generalizing the edge-level-contrast statistic would be to focus on pairs of neighbors-of-neighbors, one focal and one auxiliary, and use as the test
statistic the average outcome for focal units with treated auxiliary  neighbors-of-neighbors minus the average outcome for focal units with 
control auxiliary neighbors-of-neighbors whose treatment varies in the restricted set. In order to define the latter condition, let
$\mmp_A$ again be the set of auxiliary units, units who are not focal and who do not have any focal neighbors, and let 
 $A_i$ be an indicator for the event that unit $i$ is an auxiliary unit. Then the edge-level-contrast statistic is:
\begin{equation}\label{mbf2} T_{\rm elc}^{\rm HO}=\frac{\sum_{i,j\neq i} F_i\cdot H_{ij}\cdot A_j\cdot W_j\cdot Y^\obs_i}
{\sum_{i,j\neq i} F_i\cdot  H_{ij}\cdot A_j\cdot W_j}
-\frac{\sum_{i,j\neq i} F_i\cdot H_{ij}\cdot A_j\cdot (1- W_j)\cdot Y_i^\obs}
{\sum_{i,j\neq i} F_i\cdot H_{ij}\cdot A_j\cdot (1-W_j)}.\end{equation}

As a practical matter, tests for higher order spillovers while allowing for first order spillovers are likely to  have less power than tests for first order spillovers.  A first reason is that generally one would expect higher order spillover effects to be small relative to direct effects and first order spillover effects. Second, in the procedure discussed here, we restrict the set of assignments $\mmw_R$ that is exploited in the calculation of the p-values 
by fixing not just the assignment for focal units, but also the assignment for all their neighbors. For a given set of focal units the test for first order spillover effects would have a much larger set of auxiliary units than the test for higher order spillover effects. To counter this, it may be important to restrict the size and characteristics of the set of focal units when 
analyzing tests for higher order spillover effects.

\subsubsection{A Score Statistic}

As an alternative to the edge-level-contrast statistic, we consider a score statistic based on a 
linear-in-means model of the type considered in Manski (1993, 2013), Goldsmith-Pinkham  and Imbens (2013) and others, and previously here in Section \ref{score}. Under the null, we model the spillovers as additive and linear in the  indicator for the own treatment and the fraction of neighbors treated:
\[ Y_i^\obs=\alpha+\tau_\direct\cdot W_i+\tau_\spill\cdot
\sum_{j=1}^N W_j\cdot \overline{G}_{ij}
+\varepsilon_i,\]
where as before, $\overline{G}_{ij}=G_{ij}/\sum_{m=1}^N G_{im}$, and zero if  individual $i$ has no neighbors.

Assuming the  assignment to treatment is completely random,
we can, given this model, estimate the parameters $\alpha$, $\tau_\direct$ and $\tau_\spill$ by least squares. 
We can then consider a more general model that allows second order effects of the treatment in addition to the first order effects captured by $\tau_\spill$:
\[ Y_i^\obs=\alpha+\tau_\direct\cdot W_i+\tau_\spill\cdot
\sum_{j=1}^N W_j\cdot \overline{G}_{ij}+\tau_{\rm second}\cdot
\sum_{j=1}^N W_j\cdot \overline{H}_{ij}
+\varepsilon_i,\]
 where
$\overline{H}_{ij}=H_{ij}/\sum_{m=1}^N H_{im}$ if   $\sum_{m=1}^N H_{im}>0$, and $\overline{H}_{ij}=0$ if  $\sum_{m=1}^N H_{im}=0$.
The score statistic for the second-order spillover effect $\tau_{\rm second}$ is then proportional to the covariance between the estimated residual from this regression and the fraction of second-order neighbors who are treated:
\begin{equation}\label{corr_high} T_{\rm score}^{\rm high}=
{\rm Cov}\left(\left.
Y^\obs_i-\hat\alpha-\hat\tau_\direct\cdot W_i-\hat\tau_\spill
\cdot\sum_{j=1}^N W_j\cdot \overline{G}_{ij}
,
\sum_{j=1}^N W_j\cdot \overline{H}_{ij}
\right| \sum_{j=1}^N H_{ij}>0
\right)
.\end{equation}
This score statistic is very similar to that in the discussion of the null hypothesis of no spillovers, with two modifications. First, the outcome is now also adjusted for the first order spillover effect, by subtracting $\hat\tau_\spill
\cdot(\sum_{j=1}^N W_j\cdot \overline{G}_{ij})$, and second, we look at the correlation of this adjusted outcome with the fraction of second order neighbors who is treated, instead of the fraction of direct neighbors who is treated.

\subsection{Choosing the Focal Subpopulation for the Null Hypothesis of No Higher Order Spillovers}

Given the structure of the artificial experiment for the null of no higher order spillovers, the key to statistical power is, in
addition to the usual requirement for a sufficient number of focal units, the presence of auxilliary units (those who are not neighbors of any focal units) who are also neighbors of neighbors of focal units.   Thus, we choose the focal subpopulation to, at least approximately, maximize the number of focal-auxiliary pairs where the auxiliary unit is a neighbor of a neighbor of the 
focal unit.

Suppose we have a focal subpopulation $\mmp_F$, now with  corresponding buffer and auxiliary subpopulations $\mmp_B$ and 
$\mmp_A$. Consider adding a currently non-focal (buffer or auxiliary) individual $i$ to the focal subpopulation, changing the focal subpopulation to $\tilde \mmp_F$ and the auxiliary subpopulation to $\tilde\mmp_A$. Then $\tilde F_j=F_j$ if $j\neq i$, and $\tilde F_i=1$, $F_i=0$.
In addition, $\tilde A_i=0$, and $\tilde A_j=A_j\cdot (1-G_{ij})$ for $j\neq i$: neighbors of $i$ are removed from the set of auxiliary units.  
The number of new edges used in the edge-level-contrast statistic as a result of the change is the number of auxiliaray units that are neighbors of neighbors of $i$:
\[ \sum_{j=1}^N \tilde A_j\cdot H_{ij}=\sum_{j=1}^N  A_j\cdot(1-G_{ij})\cdot H_{ij}=\sum_{j=1}^N A_j\cdot H_{ij}.
\]
The number of old edges no longer used in the statistic after adding unit $i$ to the focal subpopulation is determined by the set of individuals who used to be auxiliary 
but become buffer units as a result of being neighbors of $i$. This leads to number of edges being dropped equal to
\[ \sum_{k=1}^N\sum_{j=1}^N  F_k\cdot (A_j-\tilde A_j)\cdot H_{kj}
+\sum_{k=1}^N F_k\cdot A_i\cdot H_{ki}\]
\[\hskip1cm
=\sum_{k=1}^N\sum_{j=1}^N  F_k\cdot  A_j\cdot G_{ij}\cdot H_{kj}
+\sum_{k=1}^N F_k\cdot A_i\cdot H_{ki}
\]
Thus, the addition of unit $i$ to the focal subpopulation 
would increase the number of comparisons by
\[ \Delta_{N,i}=\sum_{j=1}^N A_j\cdot H_{ij}
-\sum_{k=1}^N\sum_{j=1}^N  F_k\cdot  A_j\cdot G_{ij}\cdot H_{kj}
-\sum_{k=1}^N F_k\cdot A_i\cdot H_{ki}
\]
\[ \hskip1cm=\sum_{j=1}^N (A_j-A_i\cdot F_j)\cdot H_{ij}
-\sum_{k=1}^N\sum_{j=1}^N  F_k\cdot  A_j\cdot G_{ij}\cdot H_{kj}
.
\]
In cases where the alternative is proportional to the share of treated neighbors-of-neighbors, one may wish to optimize by choosing as the next focal unit the unit $i$ with the highest value for
\[ \delta_{N,i}
=\frac{\sum_{j=1}^N (A_j-A_i\cdot F_j)\cdot H_{ij}
}{\sum_{j=1}^N H_{ij}}
-\sum_{k=1}^N F_k\cdot
\frac{\sum_{j=1}^N   A_j\cdot G_{ij}\cdot H_{kj}}{\sum_{j=1}^N H_{kj}}
,\]
with the stopping rule based on whether the maximum value of $\delta_{N,i}$ over all remaining non-focal units $i$ is positive or not.

This algorithm will lead to a focal subpopulation with a large number of neighbors-of-neighbors who are auxiliary units.

\section{Exact P-values for the Null Hypothesis on Competing Network Specifications}\label{sparse}

In this section we consider null hypothesis regarding competing specifications of the network. We have two specifications of the network, $\ba_1$ and $\ba_2$, with for some pairs $(i,j)$, $G_{1,i,j}\neq G_{2,i,j}$.
We test Hypothesis \ref{sparsehyp} that
 $Y_i(\bww)=Y_i(\bww')$ for all $i$, and for  
 all pairs of assignment vectors $\bww,\bww'\in\mmw$ such that
$w_j=w'_j$ for all units $j$ such that $G_{1,ij}=1$. 

Given a set of focal units, the buffer subpopulation is now the subpopulation of units that are not focal, but that are neighbors with focal units under network $\ba_1$. The set of auxiliary units is the set of non-focal and non-buffer units.
\[ \mmv(i,\bww,H_0)=\{\bww'\in\mmw|\bww'_i=\bww_i\wedge \left(\bww'_j=\bww_j{\rm \ for\ all}\ j \ {\rm s.t.}\ G_{ij}=1\right)\}.\]
Then, as before, the restricted set of assignments is the intersection over all focal units of these sets:
\[ \mmw_S=\cap_{i\in\mmp_F} \mmv(i,\bw,H_0).\]

Next, we consider the choice of test statistics. First we consider an edge-level-contrast statistic. For all pairs of focal units and treated auxiliary units who are neighbors according to the second network, $\ba_2$, we average the outcome of the focal unit, and subtract the average, over all
all pairs of focal units and control auxiliary units who are neighbors according to the second network:
\begin{equation}\label{mbf3} T_{\rm elc}^{\rm CN}=\frac{\sum_{i,j} F_i\cdot G_{2,ij}\cdot A_j\cdot W_j\cdot Y^\obs_i}
{\sum_{i,j} F_i\cdot  G_{2,ij}\cdot A_j\cdot W_j}
-\frac{\sum_{i,j} F_i\cdot G_{2,ij}\cdot A_j\cdot (1- W_j)\cdot Y_i^\obs}
{\sum_{i,j} F_i\cdot G_{2,ij}\cdot A_j\cdot (1-W_j)}.\end{equation}
For the score statistic we first estimate the effect of spillovers from the first network as in the previous section. For focal units we then calculate the covariance of the residual from this regression with the fraction of neighbors from the second network who are treated:
\begin{equation}\label{corr_competing} T_{\rm score}^{\rm CN}=
{\rm Cov}\left(\left.
Y^\obs_i-\hat\alpha-\hat\tau_\direct\cdot W_i-\hat\tau_\spill
\cdot\frac{\sum_{j=1}^N W_j\cdot G_{1,ij}}{\sum_{j=1}^N G_{1,ij}},
\frac{\sum_{j=1}^N W_j\cdot G_{2,ij}}{\sum_{j=1}^N G_{2,ij}}
\right| \sum_{j=1}^N G_{2,ij}>0
\right)
.\end{equation}

To choose the focal subpopulation we again use a greedy algorithm, starting with the empty set at the subpopulation of focal units. We then sequentially add new focal units, one at a time, by choosing the currently non-focal unit whose inclusion in the focal subpopulation would add the most paths between focal and auxiliary units of length two, but not of length one.

\section{Exact P-values for the Null Hypothesis on Peer Effect Heterogeneity}\label{heterogeneity}

In this section we consider a null hypothesis for heterogeneity in the treatment effects,
Hypothesis \ref{noheterohyp}:
 $Y_i(\bww)=Y_i(\bww')$ for 
all  $i$, and for all pairs of assignment vectors $\bww,\bww'\in\mmw$ such that
$\sum_{j=1}^N \bww_j\cdot G_{ij}=\sum_{j=1}^N\bww'_j\cdot G_{ij}$. What we are interested in here is testing whether it matters which of one's neighbors are treated, given the number of treated neighbors. It may be that neighbors with particular characteristics are more influential than others. This maybe correspond to neighbors with similar characteristics as the ego, or neighbors who have a more central place in the network, neighbors with whom the eog has more interactions, or neighbors with particularly high values for particular characteristics.

Given a focal subpopulation, the level set is
\[ \mmv(i,\bww,H_0)=\left\{\bww'\in\mmw\left|\bww'_i=\bww_i\wedge \left(\sum_{j=1}^N\bww'_j\cdot G_{ij}=\sum_{j=1}^N\bww_j\cdot  G_{ij}\right)\right.\right\}.\]
As usual, the restricted set of assignments is the intersection over all focal units of these sets:
\[ \mmw_S=\cap_{i\in\mmp_F} \mmv(i,\bw,H_0).\]

To choose a test statistic we focus on the score approach. Under the null hypothesis we can estimate the direct and spillover effects by least squares, and calculate the residual
\[ Y^\obs_i-\hat\alpha-\hat\tau_\direct\cdot W_i-\hat\tau_\spill\cdot \frac{\sum_{j=1}^N W_j\cdot G_{ij}}{\sum_{j=1}^N G_{ij}}.\]
There is a variety of alternative hypotheses we can consider. Here we focus on one where the effect of neighbor $j$  being treated on  the outcome of individual $i$ is proportional to the degree of that unit (i.e., the number of neighbors $K_j$ that this neighbor $j$ has). This leads to
\begin{equation}\label{wat}{\rm Cov}\left(\left.
Y^\obs_i-\hat\alpha-\hat\tau_\direct\cdot W_i-\hat\tau_\spill
\cdot\frac{\sum_{j=1}^N W_j\cdot G_{ij}}{\sum_{j=1}^N G_{ij}},
\frac{\sum_{j=1}^N W_j\cdot K_j\cdot G_{ij}}{\sum_{j=1}^N K_j\cdot G_{ij}}
\right| \sum_{j=1}^N K_j\cdot G_{ij}>0
\right)
.\end{equation}

To implement this test we also need to choose the focal subpopulation. In this case it is important for focal units to have variation their friends' degree. Thus we need focal units with at least two neighbors. For each unit $i$ we calculate for all their non-focal neighbors $j$ how many non-focal neighbors this neighbor $j$ has:
\[ U_{ij}={\bf 1}_{G_{ij}=1}\cdot (1-F_j)\cdot  \sum_{j'=1, j'\neq i}^N (1-F_{j'})\cdot G_{jj'}.\]
Then we calculate the average and the standard deviation of this measure over all the neighbors of unit $i$:
\[ \overline{U}_i=\frac{\sum_{j:G_{ij=1}}  (1-F_j)\cdot U_{ij}}
{\sum_{j:G_{ij=1}}  (1-F_j)},\hskip1cm
S_{U,i}=\left(\frac{1}{K_i-1}\sum_{j:G_{ij=1}} (1-F_j) \left(U_{ij}-\overline{U}_i\right)^2
\right)^{1/2}.\]
Our approach now is to select, sequentially, focal units with high values for $S_{U,i}$.

\section{Simulations}
\label{simulations}

In this section, we carry out two sets of Monte Carlo simulations to assess the properties of the proposed procedures. In the first set, we focus on testing the null hypothesis of no spillovers in the context of general networks.
In the second, we focus on the comparison of two networks, one sparser than the other, and test the null hypothesis that all spillovers are first-order spillovers in the sparser network.

\subsection{Monte Carlo Set Up I: Testing for the Presence of Spillovers}
\label{simulations1}

The following components of the simulations are common to all designs in the first Monte Carlo set up. First consider the potential outcomes.
Let $\bww_0$ be the $N$-component vector with all elements equal to zero. Then, the baseline potential outcomes with no units exposed to the treatment are drawn from a Gaussian distribution:
\[ Y_{i}(\bww_0)\sim {\cal N}(0,1), \hskip1cm {\rm independent\ across\ all\ units}.\]
Let $\bww_{(0,i)}$ be the $N$-component vector with all elements equal to zero other than the $i$th element, which is equal to one.
We assume a constant additive direct (own) treatment effect:
\[ Y_i(\bww_{(0,i)})-Y_i(\bww_0)=\tau_{\direct},\]
for all $i=1,\ldots,N$.
Let $K_{i}$ be the number of peers for unit $i$ and let $K_{i,0}$ and $K_{i,1}$ be  the number of control and treated peers.
Then we assume a constant additive spillover effect that is proportional to the number of treated peers:
\[ Y_i(\bww)=Y_i(\bww_0)+w_i\cdot \tau_{\direct}+\frac{K_{i,1}}{K_i}\cdot \tau_{\spill}.\] If $\tau_{\spill}$ is equal to zero the null hypothesis of no spillover effects holds. If $\tau_\spill\neq 0$, the null hypothesis  is violated.

The
assignment to treatment is completely random with a fixed number of treated and control individuals. In all simulations there are 599 individuals, 300 treated individuals and 299 control individuals.

The Monte Carlo designs vary along five dimensions.
\begin{description}
  \item {\sc 1. Network Structure:}
We consider two network structures. 

In the first network structure we take a network of friendships from one of the high schools represented in the Add Health data. For details on the design of this data set see http://www.cpc.unc.edu/projects/addhealth/. We use a subset containing information on 599  students with at least one friend in the school. On average each student has 5.1 friends, with a standard deviation of 3.1, and the number of friends ranging from 1 to 18. In these simulations we keep the network fixed across the simulations. 

In the second network structure we sample Watts--Strogatz (1998) small world networks with $k = 10$ and probability of rewiring $p = 0.1$. The degree distribution thus has mean 10 and standard deviation 1.37. The size of the network is the same as in the Add Health network, 599.

\item{\sc 2. Statistic:} We consider three statistics.

The first is the edge-level-contrast statistic $T_\elc$, equal to the difference in average ego outcomes over all edges with focal egos and treated alters and the average of ego outcomes over all edges with focal egos and control alters, as given in (\ref{mbf}).
 The second is the score statistic $T_{\rm score}$ given in (\ref{corr}), motivated by a Manski-style linear-in-means model with endogenous peer effects. The third is the Aronow statistic $T_{\rm htn}$, which is the difference in average outcomes for focal units with at least one treated neighbor and those with only control neighbors.
\item{\sc 3. Own Treatment Effect:} We allow the own treatment effect $\tau_\direct$ to take on the values $0$ and $4$.
\item {\sc 4. Spillover Effect:} We allow the spillover effect $\tau_\spill$ to take on the values $0$ and $0.4$ to assess size properties of the test under the null hypothesis as well as  power of the test under the alternative hypothesis.
\item {\sc 5. Location and Number of Focal Units:} We compare three methods for choosing the focal units. In the first we randomly select 300 (approximately half) the individuals to be  focal. In the second we use the $\varepsilon$-net approach. 
In the Add Health network this approach leads to 213 (36\%) focal individuals, and in the small world networks it leads on average to   98 (16\%) focal individuals. 
In the third we maximize the number of edge comparisons, weighted by the number of neighbors, using the procedure described in Section \ref{imbensnet}. In the Add Health network this approach leads to 237 (40\%) focal individuals, and in the small world networks it leads on average to 128 (21\%) focal individuals. 
\end{description}

We approximate the p-value by drawing from the randomization distribution of the statistic under the null $1,000$ times, and calculating the proportion of of the draws where the absolute value of the statistic is larger than the absolute value of the statistic calculated on the actual data. We then report the  fraction of replications, over $4,000$ replications, where the p-value is less then 0.05. 

The results are presented in Table \ref{tab:een}. We note a couple of the findings. First of all,
when the null hypothesis is true, the tests all perform as expected, with the p-values less than 0.05 the appropriate number of times. When the null hypothesis is false we do see that the tests have substantial power. 
As discussed in the theoretical sections, the choice of focal units matters substantially for the power of the tests. Random selection of focal units performs quite poorly compared to more systematic ways of choosing the focal units. Both the method based on optimizing the number of focal-nonfocal friendships and the $\varepsilon$-net approach work substantially better.
The choice of test statistic also matters a great deal. the score statistic, designed to be optimal for interesting alternatives performs better than the edge-level-contrast statistic or the Aronow statistic. The structure of the network appears to matter less. Results for the Add Health network and the small world network are similar.

\subsection{Monte Carlo Set Up II: Testing for Sparsification}
\label{simulations2}

In the second set of simulations we focus on tests for the presence of spillovers beyond the first order spillover of a sparser network. In the simulations we take the original Add Health network with 599 students as the baseline network. We create a sparser network by randomly cutting each edge in the Add Health network with probability $q$, where either $q=0.9$ or $q=0.5$. This leads to a network with average degree 0.43 (if we cut 90\%) or   2.57 (if we cut 50\%), compared to
    5.15 in the original network.

We randomly assign 300 of the students to the treatment. We then simulate outcome data according to the linear in means model:
\[ Y_i^\obs=\tau_\direct\cdot W_i+\tau_\spill\cdot \overline{W}_{(i)}
+\varepsilon_i,\]
where $ \overline{W}_{(i)}$ is the fraction of neighbors who are treated, with weight $0\leq \lambda\leq 1$ for edges that are only present in the second, less sparse, network:
\[ \overline{W}_{(i)}=\left(\frac{\sum_{j=1}^N \left(G_{1,ij}+\lambda\cdot(G_{2,ij}-G_{1,ij})\right) W_j}{\sum_{j=1}^N \left(G_{1,ij}+\lambda\cdot(G_{2,ij}-G_{1,ij})\right)}\right).\]
If $\lambda=0$ the sparsification is appropriate because the edges only in the second network do not matter. If $\lambda=1$, the edges in the second network are just as important as those in the first network.
We simulate the $\varepsilon_i$ as independent and identically distributed, with ${\cal N}(0,1)$ distributions.

We focus on two statistics.
For the first statistic, $T_{\corr}$, based on the covariance of the residual based on the model under the null and the share of treated second-network neighbors in (\ref{corr_competing}). 
The specific statistic we focus on is the correlation between this residual and the fraction of treated neighbors for the focal individuals. The second statistic, $T_\elc$,  is the difference of two averages over all edges between focal and auxiliary individuals in (\ref{mbf3}).  The focal subpopulation is selected using the greedy algorithm described in Section \ref{sparse}.

We present results for a number of designs in Table \ref{tab:twee}. Again the test work as expected when the null hypothesis is true. The power of the test is generally higher if the spillover effect is larger $\tau_\spill=0.4$ rather than $\tau_\spill=0.1$), not surprisingly given that under the alternative the spillover effect for the second network neighbors is proportional to that for the first network neighbors. It is also higher if the sparsification of the network is more substantial ($q=0.9$ rather than $q=0.5$). Finally,  as expected  the score based statistic has more power than the edge-level-contrast.


\section{Conclusion}\label{conclusion}

In this paper we develop new methods for testing  hypotheses with experimental data in settings with a single network. We focus on the calculation of Fisher-type, exact, finite sample, p-values. The complication is that the hypotheses we are interested in are not sharp, so that conventional methods for calculating exact p-values need to be modified. We show that by analyzing an artificial experiment, different from the one actually performed, one can calculate exact p-values for interesting hypotheses regarding spillovers, sparsification of networks, and peer effect heterogeneity.  We illustrate approaches for selecting test statistics as well as the details of the artificial experiment to maximize statistical power. 
We illustrate the new methods by carrying out simulations.

\newpage

\begin{appendix}

\centerline{Appendix A: Why the Bond et al Randomization P-values are Not Valid}

\vskip0.5cm 

Bond,  Fariss,  Jones,  Kramer,  Marlow,  Settle,
and Fowler (2013), Bond et al. from hereon, are also interested in testing for spillovers (Hypothesis \ref{nopeer}). They wish to use testing procedures that are robust to the network structure.
We show here analytically that there procedures are not valid in general, and can lead to over-rejections of 0.05-level tests at rates as high as 0.20 because they ignore the variation arising from own treatment effects.

Bond et al. focus on  the  difference between the average of an ego's outcome over all edges where the alter is exposed, and the average over all edges where the alter is not exposed:
\begin{equation} T_{\bondfowler}(\bw,\by,\ba)=\frac{\sum_{i,j\neq i} G_{ij}\cdot W_j\cdot Y_i}
{\sum_{i,j\neq i} G_{ij}\cdot W_j}
-\frac{\sum_{i,j\neq i} G_{ij}\cdot (1-W_j)\cdot Y_i}
{\sum_{i,j\neq i} G_{ij}\cdot (1-W_j)}.\end{equation}
Under Hypothesis \ref{nopeer} the expected value of this statistic is zero, which makes it promising for testing this hypothesis. However,
because of the network structure there may dependence between the terms in each of these averages, and its variance is difficult to estimate for a general network structure.

Bond et al.  look at a randomization-based distribution for this statistic to test the null hypothesis of no spillovers.
The distribution is obtained by re-assigning the treatment vector $\bw$, assuming there is no effect of the treatment whatsoever, and deriving from there the quantiles of the $T_{\bondfowler}$ distribution. 
This implicitly assumes for these calculations that there is no effect of the treatment whatsoever (Hypothesis \ref{noeffects}),
which is stronger than the no-spillover null hypothesis (Hypothesis \ref{nopeer}) that they are interested in testing. The reason for this is that if one allows for direct effects of the treatment on the own outcomes, and only assumes no spillovers, one cannot infer the value of the statistic $T_{\bondfowler}$ for alternative values of the treatment assignment vector: the no-spillover null hypothesis is not sharp. The concern is that using the randomization that is based on a stronger null hypothesis is not innocuous. Bond et al justify the use of this method using simulations in which the stronger null is true.

Here we show through analytic calculations for a particular example that p-values based on these calculations are not valid, even in large samples, let alone in finite samples, and that the deviations from nominal rejection probabilities can be substantial. In general, because their calculations ignore one source of variation in the distribution of the statistic, the p-values will be too small, leading to rejections of 0.05-level tests at rates as high as 0.20.

We focus on an example with a particular network structure that allows us to simplify the large sample approximations. The population consists of $2\cdot N$ units, partitioned into $N$ pairs. Out of these $2\cdot N$ units $N$ units are  randomly selected to be exposed to the active treatment. We maintain the assumption that there are no spillovers. The potential outcomes are
\[Y_i(0)=0, \hskip1cm {\rm and}\  \ 
Y_i(1)=1, \]
so that the direct treatment effect is equal to 1.
The $N$ pairs can be partitioned into three sets: $M_{00}$ pairs with both units exposed to the control treatment, $M_{01}$ pairs with exactly one unit exposed to the control treatment and one unit exposed to the active treatment, and $M_{11}$ pairs with both units exposed to the active treatment. The number of each of these sets, $M_{00}$, $M_{01}$, and $M_{11}$ are random, but, because the total number of pairs is fixed at $N$, it follows that  $M_{00}+M_{01}+M_{11}=N$, and because exactly $N$ units are exposed to the active treatment, it must be the case that $M_{00}=M_{11}$. Hence we can rewrite these numbers in terms of a scalar random integer: define $M=M_{11}$, so that $M_{00}=M$, and $M_{01}=N-2\cdot M$.
The expected value of $M$ is $N\cdot (1/2)\cdot((N-1)/(2\cdot N))\approx N/4$. However, the variance is not $N\cdot(1/4)\cdot(3/4)$ because of the fixed number of treated units. We can approximate the large sample distribution of $\sqrt N(M/N-1/4)$ by looking at the joint distribution for 
$(\sqrt N\cdot (M_{00}/N-1/4),
\sqrt N\cdot (M_{01}/N-1/2),
\sqrt N\cdot (M_{11}/N-1/4)),$
based on independent random assignment to the treatment for each unit. This leads to
\[ \left(\begin{array}{c}
\sqrt N\cdot (M_{00}/N-1/4)\\
\sqrt N\cdot (M_{01}/N-1/2)\\
\sqrt N\cdot (M_{11}/N-1/4)\end{array}\right)\stackrel{d}{\longrightarrow} {\cal N}
\left(\left(\begin{array}{c}
0\\
0\\
0\end{array}\right)
,\frac{1}{16}\cdot \left(\begin{array}{ccc}
3 & -2& -1\\
&4 & -2\\
& & 3\end{array}\right)
\right).
\]
This implies that
\[ \left(\begin{array}{c}
\sqrt N\cdot (M_{11}/N-1/4)\\
\sqrt N\cdot (2\cdot M_{11}/N+M_{01}/N)\\
\end{array}\right)\stackrel{d}{\longrightarrow} {\cal N}
\left(\left(\begin{array}{c}
0\\
0
\end{array}\right)
,\frac{1}{16}\cdot \left(\begin{array}{cc}
3 & 4\\
 & 12\end{array}\right)
\right).
\]
Now define $M=M_{11}$ and condition on $M_{01}/N+2\cdot M_{11}/N=0$. Because the correlation between $\sqrt N\cdot (M_{11}/N-1/4)$ and $\sqrt{N}\cdot (M_{01}/N+2\cdot M_{11}/N$ is $\rho=4/sqrt{24}$, the conditional variance of 
$\sqrt N\cdot (M_{11}/N-1/4)$ given $\sqrt{N}\cdot (M_{01}/N+2\cdot M_{11}/N=0$ is $(3/16)\cdot (1-\rho^2)=1/16$,
and 
\[ \sqrt N\cdot \left(\frac{M}{N}-\frac{1}{4}\right)\stackrel{d}{\longrightarrow}
{\cal N}\left(0,\frac{1}{16}\right).\]

Now consider the statistic $T_{\bondfowler}$. We calculate first the actual distribution of this statistic under the randomization distribution. Then we compare this to the distribution Bond et al  use for the calculation of p-values.

There are $2\cdot N$ edges. Out of these $N$ have treated alters and $N$ have control alters. For the $N$ edges with treated alters $2\cdot M_{11}=2\cdot M$ have treated egos, and so have realized outcome equal to $Y_i(1)=1$, and $M_{01}=N-2\cdot M$ have control egos, and so have realized outcomes equal to $Y_i(0)=0$. The average realized outcome for egos with treated alters is therefore $2\cdot M/N$.
Similarly, for the $N$ edges with control alters, there are $2\cdot M_{00}=2\cdot M$  edges with control egos and realized outcomes $Y_i(0)=0$, and $M_{01}=N-2\cdot M$ edges with treated egos and thus $Y_i(1)=1$, leading to an average realized outcome equal to $1-2\cdot M/N$. Hence the value of the statistic is
\[ T_{\bondfowler}=2\cdot \frac{M}{N}-\left(1-2\cdot\frac{M}{N}\right)=4\cdot\left(\frac{M}{N}-\frac{1}{4}\right).\] 
The actual distribution of the normalized statistic, under random assignment, is
\[ \sqrt{N}\cdot T_{\bondfowler}
=\sqrt{N}\cdot \left(\frac{4\cdot M}{N}-1\right)
\stackrel{d}{\longrightarrow}
 {\cal N}\left(0,1\right)
.\]

Now consider the distribution used by Bond et al for the calculation of their p-values. They calculate the randomization distribution, assuming that there are no effects of the treatment whatsoever. 
Under this randomization distribution, there are always $N$ egos with treated alters, and $N$ egos with control alters. Out of the $2\cdot N$ units there are $N$ with realized outcome equal to 1 and $N$ with realized outcome equal to 0, so that the total average outcome is exactly $1/2$. Hence, if the average of the outcome for the egos with treated alters is equal to $\overline{Y}_{\rm t}$, the average of the outcome for egos with control alters is equal to $\overline{Y}_{\rm c}=1-\overline{Y}_{\rm t}$. Therefore the difference in the average outcome for egos with treated alters and the average outcome for egos with control alters is equal to $2\cdot\overline{Y}_{\rm t}-1$. 
To infer the randomization distribution used by Bond et al, we need to infer the distribution of $\overline{Y}_{\rm t}$ under their randomization distribution. 
We can write $\overline{Y}_{\rm t}$ as
\[ \overline{Y}_{\rm t}=\frac{1}{N}\sum_{i=1}^{2N} W^p_i\cdot Y_i,\]
where $W^p_i$ is an indicator for unit $i$ having a treated alter. We are interested in this distribution under random assignment of $Z_i$, with $\sum_{i=1}^{2N} Z_i=N$, for fixed ${\bf Y}$.
(It is the treating of $\by$ as fixed that is not correct here -- if we change the treatment of the alter  for unit $i$ we may be changing the value of the outcome for unit$i$'s alter. Thus the $Y_i$ are stochastic, leading to additional variation in the test statistic that is not taken into account in the B procedure.) Note that $\sum_{i=1}^{2\cdot N} Y_i=N$ and $\sum_{i=1}^{2\cdot N} W^p_i=N$. The treatments (and thus the peer treatments) are randomly assigned, with ${\rm pr}(W^p_i=1)=1/2$ and ${\rm pr}(W^p_i=1|W^p_j=1)=(N-1)/(2\cdot N-1)$. Define $D_i=2\cdot W^p_i-1$ so that $W_i^p=(D_i+1)/2$, and
\[ \mme[D_i]=0,\hskip1cm D^2_i=1,\hskip1cm \mme[D_i\cdot D_j]=-\frac{1}{2\cdot N-1},\ {\rm for}\ j\neq i.\]
Now
\[ \overline{Y}_{\rm t}=\frac{1}{N}\sum_{i=1}^{2N}  Y_i\cdot \frac{D_i+1}{2}\]
\[\hskip1cm =\frac{1}{N}\sum_{i=1}^{2N}  Y_i\cdot \frac{1}{2}
+\frac{1}{2N}\sum_{i=1}^{2N}  Y_i\cdot D_i=\frac{1}{2}+\frac{1}{2N}\sum_{i=1}^{2N}  Y_i\cdot D_i.\]
Then
\[ \mme\left[\overline{Y}_{\rm t}\right]=1/2,\]
and
\[ \mmv\left(\overline{Y}_{\rm t}\right)=\frac{1}{4\cdot N^2}\cdot\mme\left[\left(\sum_{i=1}^{2N} Y_i\cdot D_i\right)^2\right]
=\frac{1}{4\cdot N^2}\cdot\mme\left[\sum_{i=1}^{2N} D_i^2\cdot Y_i^2+
\sum_{i=1}^{2N}\sum_{j\neq i} D_i\cdot D_j\cdot Y_i\cdot Y_j
\right]
\]
\[\hskip1cm =\frac{1}{4\cdot N^2}\cdot \sum_{i=1}^{2N}  Y_i+
\frac{1}{4\cdot N^2}\cdot\sum_{i=1}^{2N}\sum_{j\neq i} Y_i\cdot Y_j\cdot \mme[D_1\cdot D_2]
\]
\[\hskip1cm =\frac{1}{4\cdot N}-
\frac{1}{4\cdot N^2}\cdot N\cdot (N-1)\cdot\frac{1}{2\cdot N-1}
\]
\[\hskip1cm =\frac{1}{4\cdot N}-
\frac{1}{4\cdot N}\cdot \frac{N-1}{2\cdot N-1}\approx \frac{1}{8\cdot N}.
\]
Hence the variance of $N\cdot \overline{Y}_{\rm t}$ is equal to $1/8$, and thus the variance of Bond et al randomization distribution is $4\cdot N\cdot \mathbb{V}( \overline{Y}_{\rm t})$ which is equal to 0.5. The actual distribution has variance equal to 1, which is twice as large. The implication is that the for a two-sided test at the 0.05 level the rejection probability based on using the incorrect Bond et al randomization distribution is 0.157. Bond et al implicitly use the wrong variance of 0.5 for the test statistic, leading to 
\[ {\rm pr}\left(\sqrt{2}\cdot |T_\bondfowler|>1.96\right)
={\rm pr}\left(|T_\bondfowler|>\sqrt{2}\cdot 1.96\right)\]
\[\hskip1cm ={\rm pr}\left(|T_\bondfowler|>\frac{ 1.96}{\sqrt{2}}\right)
\approx{\rm pr}\left(|T_\bondfowler|>1.386\right)\approx 0.157.
\]

We carried out a small simulation study to verify these analytic calculations. We use $N=1000$ pairs,  10,000 replications, and use 1,000 draws from the randomization distribution. We reject the null hypothesis if the  Bond et al p-value is less than 0.05. This leads us to reject at a rate equal to 0.153, close to the theoretical rejection rate we calculated above which is equal to 0.157. (A 95\% confidence interval for the rejection rate is $(0.144, 0.163)$).

\end{appendix}

\renewcommand{\appendix}{\footnotesize\parindent 0cm\setcounter{equation}{0} 
\renewcommand{\theequation}{B.\arabic{equation}}
\setcounter{lemma}{0}\renewcommand{\thelemma}{B.\arabic{lemma}}}

\vskip0.5cm 
\begin{appendix}

\centerline{Appendix B: Derivation of the Score Test Statistic for the Null of No Spillovers}

\vskip0.5cm 
In terms of the potential outcomes the linear-in-means model in (\ref{lim}) 
 corresponds to
\begin{equation} \by(\bww)=
\alpha_0\cdot \left(I-\tau_{\rm endog}\cdot\overline{\ba}\right)^{-1}\cdot\iota_N+
\tau_\direct\cdot \left(I-\tau_{\rm endog}\cdot\overline{\ba}\right)^{-1} \bww+\left(I-\tau_{\rm endog}\cdot\overline{\ba}\right)^{-1} \varepsilon.
\end{equation}
The expected value of the observed outcomes given the assignment is, given the random assignment,
\begin{equation}\label{expectation} \mme[\by^\obs|\bw=\bww]=
\mme[\by(\bww)]=
\alpha_0\cdot \left(I-\tau_{\rm endog}\cdot\overline{\ba}\right)^{-1}\iota_N+
\tau_\direct\cdot \left(I-\tau_{\rm endog}\cdot\overline{\ba}\right)^{-1} \bww.\end{equation}
Under the null hypothesis that $\tau_{\rm endog}=0$, the least squares estimates for the remaining parameters based on outcomes for focal units are 
\[ \hat\alpha_0=\overline{Y}^\obs_{F,0},\hskip1cm {\rm and}
\ \ \hat\tau_\direct=\overline{Y}^\obs_{F,1}
-\overline{Y}^\obs_{F,c},\]
where, for $w=0,1$, $\overline{Y}^\obs_{F,w}$ is the average outcome for  focal units with $W_i=w$,
\[\overline{Y}^\obs_{F,w}=\frac{1}{N_{F,w}}\sum_{i:F_i=1,W_i=w}  Y^\obs_i,,
\]
and $N_{F,w}$ is the number of focal units with $W_i=w$.
Hence the residual under the null is
\[ \hat\varepsilon_i^{\rm null}=Y^\obs_i-\hat\alpha_0-W_i\cdot\hat\tau_\direct.\]
Under normality of the outcome the score for $\tau_{\rm endog}=0$ is proportional to the covariance of the residual under the null and the
 derivative of the expectation in (\ref{expectation}), with respect to $\tau_{\rm endog}$, evaluated at $\tau_{\rm endog}=0$. The derivative of the expectation at $\tau_{\rm endog}=0$ is
\[ \frac{\partial}{\partial\tau_\direct}\mme[\by^\obs|\bw]=
\alpha_0\cdot\overline{\ba}
\iota_N
+
\tau_\direct\cdot \overline{\ba}\bw=
\alpha_0\cdot \overline{\ba}(\iota_N-\bw)+
(\tau_\direct+\alpha_0)\cdot\overline{\ba}\bw
.\]
Substituting $\overline{Y}^\obs_{F,0}$ for $\alpha_0$ and 
$\overline{Y}^\obs_{F,1}-\overline{Y}^\obs_{F,0}$ for $\tau_\direct$ suggests that
 a natural test statistic would be the covariance  of the residual under the null and 
$\overline{Y}^\obs_{F,0}\cdot \overline{\ba}(\iota_N-\bw)+\overline{Y}^\obs_{F,1}\cdot \overline{\ba}\bw$.
This leads to the following average score:
\[ \frac{1}{N_F}\sum_{i\in\mmp_F} \left\{\left( Y^\obs_i-\overline{Y}^\obs_{F,0}
-W_i\cdot (\overline{Y}^\obs_{F,1}-\overline{Y}^\obs_{F,0})
\right)\cdot 
\sum_{j=1}^N \overline{G}_{ij}\cdot \Bigl(
(1-W_j)\cdot \overline{Y}^\obs_{F,0}+W_j\cdot \overline{Y}^\obs_{F,1}\Bigr)
\right\}
.\]
Because $\sum_{j=1}^N \overline{G}_{ij}=1$, in combination with the fact that the residuals average to zero, it follows that the score statistic is proportional to the covariance between the residual under the null and $\sum_{j=1}^N \overline{G}_{ij}\cdot W_j$, which is the fraction of treated neighbors, leading to the score statistic
\[ T_{\rm score}=
{\rm Cov}\left(\left.Y^\obs_i-\overline{Y}^\obs_{F,0}
-W_i\cdot (\overline{Y}^\obs_{F,1}-\overline{Y}^\obs_{F,0}),
\sum_{j=1}^N W_j\cdot \overline{G}_{ij}
\right| \sum_{j=1}^N G_{ij}>0,F_i=1
\right)\]
\[\hskip2cm
=
{\rm Cov}\left(\left.
Y^\obs_i-\hat\alpha-\hat\tau_\direct\cdot W_i,
\sum_{j=1}^N W_j\cdot \overline{G}_{ij}
\right| \sum_{j=1}^N G_{ij}>0,F_i=1
\right)
,\]
which is the expression in (\ref{corr}).
\end{appendix}

\newpage

\centerline{\Large{\bf References}}
\vskip0.5cm

\begin{description}


\item[]{\small \textsc{Aronow, P.,} (2012),
``A general method for detecting interference between units in randomized experiments,''
{\it 
Sociological Methods \& Research}, Vol. 41(1): 3-16.}

\item[]{\small \textsc{Aronow, P., and C. Samii}, (2013),
``Estimating Average Causal Effects Under Interference Between Units'', unpublished working paper.}

\item[]{\small \textsc{Baird, S., A. Bohren, C. McIntosh, and B. Ozler,} (2014),
``Designing Experiments to Measure Spillover Effects,'' Policy Research Working Paper 6824, The World Bank.}

\item[]{\small \textsc{Barrios, T., R. Diamond, G. Imbens, and M. Koles\'ar} (2012),
``Clustering, Spatial Correlations and Randomization Inference'' {\it Journal of the American Statistical Association}, Vol. 107(498):578-591.}

\item[]{\small \textsc{Basu, D.,} (1980),
``Randomization Analysis of Experimental Data: The Fisher Randomization Test,''  {\it Journal of the American Statistical Association}, Vol. 75(371): .}

\item[]{\small \textsc{Bond, R., C. Fariss, J. Jones, A. Kramer, C. Marlow, J. Settle
and J. Fowler} (2012),
``A 61-million-person experiment in social influence
and political mobilization'' {\it Nature}, 295-298.}

\item[]{\small \textsc{Bowers, J., M. Fredrickson, and C. Panagopoulos,} (2012),
``Reasoning about Interference Between Units: A General Framework,'' {\it Political Analysis}, 21:97-124.}

\item[]{\small \textsc{Canay, I., J. Romano, and A. Shaikh,} (2015),
``Randomization Tests under an Approximate Symmetry Assumption,''
Unpublished Manuscript.}

\item[]{\small \textsc{Carrell, S.,
B. Sacerdote, and
J. West} (2013),
``From Natural Variation to Optimal Policy? The Importance of Endogenous Peer Group Formation'' {\it Econometrica}, 81(3): 855-882.}

\item[]{\small \textsc{Choi, D.}, (2014),
``Estimation of Monotone Treatment Effects in Network Experiments,''  arXiv:1408.4102v2.}

\item[]{\small \textsc{Christakis, N., and Fowler}, (2007),
``The Spread of Obesity in a Large Social Network over 32 Years,''  {\it The New England Journal of  Medicine}, 357:370-379.}

\item[]{\small \textsc{Eckles, D., B. Karrer, and J. Ugander} (2014),
``Design and Analysis of Experiments in Networks: Reducing Bias from Interference'', unpublished working paper.}

\item {\small \textsc{Fisher, R. A.}, (1925), \emph{Statistical Methods for Research Workers}, 1st ed, Oliver and Boyd, London.}

\item {\small \textsc{Fisher, R. A.}, (1935), {\it Design of Experiments}, Oliver and Boyd. }

\item[]{\small \textsc{Goldenberg, A., A. Zheng, S.
Fienberg and E. Airoldi,} (2009),
``A Survey of Statistical Network Models,'' {\it Foundations and Trends in
Machine Learning}
Vol. 2(2): 129�233.}

\item[]{\small \textsc{Goldsmith-Pinkham, P., and G. Imbens,} (2013),
``Social networks and the identification of peer effects'' (with discussion), {\it Journal of Business and Economic Statistics}, 31(3):253-264.}

\item[]{\small \textsc{Gupta, A., R. Krauthgamer, and J R. Lee,} (2003) ``Bounded geometries, fractals, and low-distortion embeddings.'' In {\it Foundations of Computer Science, 2003. Proceedings. 44th Annual IEEE Symposium on}, pp. 534-543. IEEE, 2003.}

\item[]{\small \textsc{Holland, P. and S. Leinhardt,} (1981),
``An Exponential Family of Distributions for Directed Graphs'' (with discussion), {\it Journal of the American Statistical Association}, 76:33-65.}

\item[]{\small \textsc{Hudgens, M., and M. Halloran,} (2008),
``Toward Causal Inference With Interference'' {\it Journal of the American Statistical Association}, 103(482): 832-842.}

\item[]{\small \textsc{Imbens, G., and D. Rubin,} (2015),
 {\em Causal Inference in Statistics, Social, and Biomedical Sciences: An Introduction}, Cambridge University Press, Cambridge, UK.}


\item[]{\small \textsc{Kolaczyk,} (2009),
 {\em Statistical Analysis of Network Data: Methods and Models}, Springer Verlag.}

\item[]{\small \textsc{Lehmann, E., and J. Romano}, (2005), \emph{Testing Statistical Hypotheses},
Springer Verlag, New York. }

\item[]{\small \textsc{Li, H., and E. Loken} (2002),
``A Unified Theory of Statistical Analysis and Inference for Variance Component Models for Dyadic Data,'' {\em Statistica Sinica}, 12, 519-535.}

\item[]{\small \textsc{Liu, L., and M. Hudgens} (2013),
``Large Sample Randomization Inference of Causal Effects in the Presence of Interference,'' {\em Journal of the American Statistical Association}, 288-301.}

\item[]{\small \textsc{Manski, C.,} (1993),
``Identification of Endogenous Social Effects: The Reflection Problem,'' {\em Review of Economic Studies}, 60(3):531-542.}

\item[]{\small \textsc{Manski, C.,} (2013),
``Identification of Treatment Response with Social Interactions,'' {\em The Econometrics Journal}, 16(1): S1-S23.}

\item[]{\small \textsc{Manski, C., and E. Tamer} (2002), ``Inference on regressions with interval data on a regressor or outcome,''
{\it Econometrica} Vol. 70:519--47.}

\item[]{\small \textsc{Newman, M. E.} (2006), ``Modularity and community structure in networks,'' {\it Proceedings
of the National Academy of Sciences}, 103, 8577--8582.}

\item[]{\small \textsc{Ogburn, E. L., and T. J. VanderWeele} (2014), ``Causal diagrams for interference,'' {\it Statistical Science} 29(4)Munmun De Choudhury: 559-578.}

\item[]{\small \textsc{Rosenbaum, P.}, (1984), ``Conditional Permutation Tests and the Propensity Score in Observational Studies'', \emph{Journal of the American Statistical Association}, Vol. 79(387): 
565-574. }

\item[]{\small \textsc{Rosenbaum, P.}, (1995, 2002), \emph{Observational Studies},
Springer Verlag, New York. }

\item[]{\small \textsc{Rosenbaum, P.}, (2007), ``Interference between units in randomized experiments,'' \emph{Journal of the American Statistical Association}, Vol.  102: 191-200.}

\item[]{\small \textsc{Rosenbaum, P.}, (2009), \emph{Design of Observational Studies}, 
Springer Verlag, New York. }

\item[]{\small\textsc{Rubin, D.} (1974), "Estimating Causal Effects of Treatments in
Randomized
and Non-randomized Studies,\textquotedblright {\em Journal of Educational Psychology}, 66,
688-701. }

\item[]{\small \textsc{Rubin, D.,} (1980),
``Randomization Analysis of Experimental Data: The Fisher Randomization Test,'' Comment {\it Journal of the American Statistical Association}, Vol. 75(371): 591-593.}

\item[]{\small \textsc{Sacerdote, B.,} (2001),
``Peer Effects with Random Assignment: results for Dartmouth Roommates,'' {\em Quarterly Journal of Economics}, 116(2):681-704.}

\item[]{\small \textsc{Tchetgen Tchetgen, E., and
T. VanderWeele,} (2012),
``On causal inference in the presence of interference,'' {\it Statistical Methods in Medical Research}, vol. 21  no. 1  55-75.}

\item[]{\small \textsc{Thomas, A., and J. Blitzstein}, (2011),
``Valued Ties Tell Fewer Lies: Why Not To Dichotomize
Network Edges With Thresholds,'' unpublished manuscript.}

\item[]{\small \textsc{Toulis, P., and E. Kao},  (2013), ``Estimation of Causal Peer Influence Effects,'' {\it JMLR W\&CP} 28(3): 1489�1497.}

\item[]{\small \textsc{Ugander, J., B. Karrer, L. Backstrom, and J. Kleinberg} (2013),
``Graph Cluster Randomization: Network Exposure to Multiple Universes,''
{\it Proc. of KDD}. ACM.}

\item[]{\small \textsc{Van Der Laan, M.,} (2014),
``Causal Inference for a Population of Causally Connected Units,'' {\em Journal of Causal Inference}, 1-61.}

\item[]{\small \textsc{Watts, D.,  and S. Strogatz,} (1998), ``Collective dynamics of 'small-world' networks,'' {\it Nature} Vol. 393(6684): 440�442.} 

\end{description}

\newpage

\begin{table}[ht]
  \centering
  \caption{\sc Rejection Rates of Null Hypothesis of No Spillovers}
  \vskip0.5cm
  \begin{tabular}[ht]{ccccccc}
\hline\\
&  &  Own & Spillover &\multicolumn{3}{c}{Focal Vertex Selection}\\
Network& Statistic & Effect &  Effect & Random & $\varepsilon$-net & $\delta_{N,i}$\\  \\ \hline \\
Add Health& $T_{\rm score}$ & 0 & 0  & 0.059& 0.056& 0.045\\ 
& $T_\elc$  & 0 & 0  & 0.058& 0.054& 0.044\\ 
&  $T_{\rm htn}$  & 0 & 0  & 0.059& 0.039 & 0.046\\  \\
&$T_{\rm score}$ & 4 & 0  & 0.056&0.053&0.051\\ 
& $T_\elc$ & 4 & 0  & 0.051&0.048 & 0.059 \\ 
& $T_{\rm htn}$  & 4 & 0  & 0.050&0.053 & 0.051\\ \\
&$T_{\rm score}$& 0 & 0.4  & 0.362& 0.463 & 0.527\\ 
& $T_\elc$ & 0 & 0.4  & 0.174&0.299&0.413\\ 
&  $T_{\rm htn}$  & 0 & 0.4  & 0.141&0.296&0.327\\ \\ 
& $T_{\rm score}$ & 4 & 0.4  & 0.346& 0.461 & 0.529\\ 
& $T_\elc$ & 4 & 0.4  & 0.083& 0.102 & 0.123 \\ 
&  $T_{\rm htn}$  & 4 & 0.4  & 0.069& 0.088 & 0.116\\ 
\\
Small World & $T_{\rm score}$ & 0 & 0  & 0.046& 0.048& 0.054 \\ 
($K=10,p_{\rm rw}=0.1$)& $T_\elc$ & 0 & 0  & 0.048&0.040& 0.057\\ 
& $T_{\rm htn}$ & 0 & 0  & 0.041& 0.049& 0.050\\ \\
& $T_{\rm score}$& 4 & 0  & 0.055& 0.046& 0.050\\ 
& $T_\elc$ & 4 & 0  & 0.049& 0.054& 0.055\\ 
& $T_{\rm htn}$  & 4 & 0  & 0.053& 0.054&0.044\\ \\
&$T_{\rm score}$ &0 & 0.4  & 0.155& 0.090& 0.131\\ 
& $T_\elc$ & 0 & 0.4  & 0.112&0.092& 0.128\\ 
&  $T_{\rm htn}$  & 0 & 0.4  & 0.059&0.042& 0.065\\  \\
& $T_{\rm score}$& 4 & 0.4  & 0.153& 0.095& 0.154\\ 
& $T_\elc$ & 4 & 0.4  & 0.060& 0.060& 0.061\\ 
&  $T_{\rm htn}$  & 4 & 0.4  & 0.047& 0.047& 0.050\\ 
\\
\hline\\
     \end{tabular}
  \label{tab:een}
\end{table}

\begin{table}[ht]
  \centering
  \caption{\sc Rejection Rates of Null Hypothesis of No Spillovers Beyond the First Order Spillovers from the Sparsified Network, AddHealth data, 10,000 Replications}
  \vskip0.5cm
  \begin{tabular}[ht]{cccccc}
\hline\\
  &   &  &&\multicolumn{2}{c}{Prop of Links Dropped}
\\
 Statistic & $\tau_\direct$ &  $\tau_\spill$ &$\lambda$& $q=0.9$& $q=0.5$
\\  \\ \hline \\
\\ 
$T_\corr$ & 0 & 0.1  & 0&  0.051  
&   0.051  
\\
 $T_{\elc}$ & 0 & 0.1  & 0& 0.050  
& 0.049 
\\
 $T_\corr$ & 0 & 0.4  & 0& 0.051 
& 0.050 
\\ 
 $T_{\elc}$ & 0 & 0.4  & 0& 0.050
&  0.050 
\\ \\
$T_\corr$ & 4 & 0.1  & 0& 0.052  
&  0.046  
\\
 $T_{\elc}$ & 4 & 0.1  & 0& 0.049  
& 0.046 
\\
 $T_\corr$ & 4 & 0.4  & 0&  0.058 
&  0.048
\\ 
 $T_{\elc}$ & 4 & 0.4  & 0& 0.051 
&  0.047 
\\ \\
 $T_\corr$ & 0 & 0.1  & 0.5& 0.060 & 0.055
\\
 $T_\elc$ & 0 & 0.1  & 0.5& 0.054 & 0.048
\\ 
 $T_\corr$ & 0 & 0.4  & 0.5& 0.212 & 0.108
\\
 $T_\elc$ & 0 & 0.4  & 0.5& 0.121 & 0.069
\\ 
\\
 $T_\corr$ & 4 & 0.1 & 0.5& 0.057 & 0.053
\\
 $T_\elc$ & 4 & 0.1 & 0.5 & 0.052 
&   0.047 
\\
 $T_\corr$ & 4 & 0.4  & 0.5& 0.212
& 0.112
\\
 $T_\elc$ & 4 & 0.4  & 0.5& 0.061 
& 0.051
\\ 
\\
\hline\\
     \end{tabular}
  \label{tab:twee}
\end{table}

\end{document}